\documentclass[reqno]{amsart}

\usepackage{mathtools}

\usepackage[demo]{graphicx}
\usepackage{caption}
\usepackage{subcaption}

\usepackage{tikz-cd}
\usetikzlibrary{positioning}

\usepackage{tcolorbox}

\usepackage{enumerate}
\usepackage{amsmath}
\usepackage{amssymb}
\usepackage{graphicx}
\usepackage[cmtip,all]{xy}
\usepackage[mathscr]{eucal}
\usepackage{bbold}            
\usepackage{tikz}[all]
\usepackage{tikz-cd}

\usepackage{url}
\usepackage{hyperref}
\hypersetup{
  colorlinks   = true, 
  urlcolor     = blue,
  linkcolor    = blue, 
  citecolor   = blue 
}

\usepackage{mathrsfs}

\newtheorem{theorem}{Theorem}[section]

\newtheorem{corollary}[theorem]{Corollary}
\newtheorem{proposition}[theorem]{Proposition}

\newtheorem{claim}{Claim}
\newtheorem*{theorem*}{Theorem}

\theoremstyle{definition}

\newtheorem{definition}[theorem]{Definition}
\newtheorem{example}[theorem]{Example}

\newtheorem{notation}[theorem]{Notation}

\newtheorem{fact}[theorem]{Fact}

\theoremstyle{remark}
\newtheorem{remark}[theorem]{Remark}

\numberwithin{equation}{section}

\newcommand{\Lwedge}{\mathsf{\Lambda}}

\newcommand{\ad}[1]{\mathsf{ad}(#1)}

\newcommand{\field}[1]{\ensuremath{\mathbb{#1}}}
\newcommand{\bba}{\field{A}}
\newcommand{\bbc}{\field{C}}

\newcommand{\bbg}{\field{G}}

\newcommand{\bbl}{\field{L}}

\newcommand{\bbn}{\field{N}}

\newcommand{\bbr}{\field{R}}

\newcommand{\bbt}{\field{T}}
\newcommand{\bbu}{\field{U}}

\newcommand{\bbx}{\field{X}}

\newcommand{\bbz}{\field{Z}}

\newcommand{\place}[1]{\ensuremath{\mathbb{#1}}}

\newcommand{\bone}{\place{1}}
\newcommand{\cala}{\mathcal{A}}
\newcommand{\calb}{\mathcal{B}}
\newcommand{\calc}{\mathcal{C}}
\newcommand{\cald}{\mathcal{D}}

\newcommand{\cali}{\mathcal{I}}
\newcommand{\calk}{\mathcal{K}}

\newcommand{\cals}{\mathcal{S}}

\newcommand{\scrA}{\mathscr{A}}
\newcommand{\scrB}{\mathscr{B}}
\newcommand{\scrC}{\mathscr{C}}
\newcommand{\scrD}{\mathscr{D}}

\newcommand{\scrG}{\mathscr{G}}

\newcommand{\scrI}{\mathscr{I}}

\newcommand{\scrO}{\mathscr{O}}

\newcommand{\scrV}{\mathscr{V}}

\newcommand{\frg}{\mathfrak{g}}

\newcommand{\frh}{\mathfrak{h}}

\newcommand{\frl}{\mathfrak{l}}

\newcommand{\frakm}{\mathfrak{m}}

\newcommand{\frs}{\mathfrak{s}}

\newcommand{\fru}{\mathfrak{u}}

 \newcommand{\mba}{\mathbf{a}}
 
 \newcommand{\mbb}{\mathbf{b}}
 
 \newcommand{\mbc}{\mathbf{c}}
 \newcommand{\mbd}{\mathbf{d}}
 \newcommand{\mbe}{\mathbf{e}}

 \newcommand{\mbm}{\mathbf{m}}
 \newcommand{\mbn}{\mathbf{n}}
 
 \newcommand{\mbp}{\mathbf{p}}

  \newcommand{\mbu}{\mathbf{u}}

 \newcommand{\mbv}{\mathbf{v}}

 \newcommand{\mbx}{\mathbf{x}}
 
\newcommand{\mby}{\mathbf{y}}

 \newcommand{\mat}[2]{\mathcal{M}_{#1 \times #2}}
 \newcommand{\la}{\langle}
 \newcommand{\ra}{\rangle}

\newcommand{\spec}[1]{\mathsf{Spec}\left(#1 \right)}

\newcommand{\sch}{\mathsf{Sch}}  
\newcommand{\Hom}[3]{\mathsf{Hom}_{#1}(#2, #3) } 
\newcommand{\HomSch}[2]{\underline{\mathsf{Hom}}(#1, #2) } 

\newcommand{\proj}[1]{\mathsf{proj_{#1}}}

\newcommand{\vtriang}{\Delta \hspace{-0.235cm}\mathord{\raisebox{-0.2\depth}{\scalebox{0.42}{\( \Delta\)} } }
\hspace{-0.26cm}\mathord{\raisebox{-0.8\depth}{\scalebox{0.64}{\( \Delta\)}}}} 

\newcommand{\simplex}[1]{\vtriang{\kern.22em}^{#1}} 
\makeatletter
\newcommand{\extp}{\small \@ifnextchar^\@extp{\@extp^{\,}}}
\def\@extp^#1{\mathop{\bigwedge\nolimits^{\!#1}}}
\makeatother

\makeatletter
\newcommand*{\Relbarfill@}{\arrowfill@\Relbar\Relbar\Relbar}
\newcommand*{\xeq}[2][]{\ext@arrow 0055\Relbarfill@{#1}{#2}}
\makeatother

\newcommand{\lulu}{\textsc{lulu}}
\newcommand{\Ad}[2]{\textsf{Ad}_{#1}(#2) }

\begin{document}

\title[LULU is syntomic]{LULU is syntomic}


\author{E. Javier Elizondo}
\address{Instituto de Matem\'aticas, Universidad Nacional Aut\'onoma de M\'exico (UNAM)}
\thanks{The first author would like to thank the fellowship granted by DGAPA of the National University of Mexico (UNAM). He is also grateful for the hospitality of the Department of Mathematics, Texas A\&{M} University, during his sabbatical year}

\author{Paulo Lima-Filho}
\address{Department of Mathematics, Texas A\&{M} University}
\thanks{}


\begin{abstract}
Let $G$ be a Chevalley group over a field $\Bbbk$. Fix a maximal torus \(\mathbb{T} \) in \( G \), along with opposite Borel subgroups \( B \) and \( B^-\) satisfying \( \bbt = B \cap B^-\), and denote by \( \bbu := R_u(B) \) and \( \bbl := R_u(B^-)\) their respective unipotent radicals. We prove that the multiplication map $ \mu \colon \bbl \times \bbu \times \bbl \times \bbu \longrightarrow G $  is syntomic and faithfully flat over any base field \( \Bbbk\).
\end{abstract}

 \maketitle

\setcounter{tocdepth}{1}
\tableofcontents

%


\section{Introduction}
\label{sec:intro}

The decomposition of a square matrix into the product of a lower triangular and an upper triangular matrix is a standard process in computational linear algebra and in various applications of linear algebra. The range of such applications vary from numerical analysis to signal processing and computer graphics; see \cite{MR3292660}, \cite{MR2132485} and \cite{MR1875822}.

In a more specific context, it has been shown in \cite{MR1466897} that any matrix \( A \in SL_n(\Bbbk) \), where \( \Bbbk \) is a field, can be expressed as a product
\( A=L_1 U_1 L_2 U_2\), where
\( L_1, L_2\) are elements of the subgroup \( \bbl \) of lower triangular unipotent matrices (i.e., lower triangular matrices with \( 1\)'s along the diagonal), and \( U_1, U_2\) belong to the subgroup \( \bbu\) of upper triangular unipotent matrices. This type of decomposition is referred to as a \lulu\ decomposition of \( A\),
and it is sometimes rephrased as ``every `unit' matrix can be written as a product of at most four generalized shears.''

To broaden this perspective, let \( \Phi\) be a reduced, irreducible root system and  \( G_\Phi \) denote the associated simply-connected Chevalley group (a split semi-simple linear algebraic group over  \( \bbz \)). Given a Borel subgroup  \( B_\Phi\) containing a maximal torus \( \bbt\), let    \( U_\Phi\) denote its  unipotent radical and \( U^-_\Phi\) denote the unipotent radical of the opposite Borel subgroup. Using techniques from Bass \cite{MR174604} and  Tavgen \cite{MR1175793}, it is shown in \cite{MR2822515} that for any commutative ring \( R \) with stable rank \( 1\), any element \( A \in G_\Phi(R) \) can be decomposed as a product \( A = L_1 U_1 L_2 U_2,\) with \( L_i \in U^-_\Phi(R)\) and \( U_i \in U_\Phi(R),\, i=1,2\). 
In particular,  it follows that the   map   
\( U^-_{\Phi, \Bbbk}\times U_{\Phi, \Bbbk}\times U^-_{\Phi, \Bbbk}\times U_{\Phi, \Bbbk}  \to G_{\Phi, \Bbbk} \) of schemes over a field \( \Bbbk\) given by 
multiplication, is surjective. Whenever \( \Phi \), \( \Bbbk \) and \( B_\Phi\) are understood, we denote 
\( \bbu := U_{\Phi, \Bbbk}, \ \bbl := L_{\Phi, \Bbbk}\) and \( G:= G_{\Phi, \Bbbk} \), and call the multiplication map
\begin{equation}
\label{eq:LULU}
\mu \ \colon \ \bbl \times \bbu \times \bbl\times \bbu  \longrightarrow G,
\end{equation}
the \lulu\ map.

Let us recall that (according to Fontaine and Messing), the term \emph{syntomic} was coined by Barry Mazur to designate a morphism of schemes which  is \emph{flat of finite presentation and a local complete intersection morphism}; see 
\cite[\href{https://stacks.math.columbia.edu/tag/01UB}{Section 01UB}]{stacks-project}, \cite[\href{https://stacks.math.columbia.edu/tag/068E}{Lemma 37.62.8}]{stacks-project}, and   also \cite[Defn. 19.3.6]{MR173675}. These morphisms are used to define the syntomic topology, and the notion of  syntomic cohomology, which plays important roles in the study of \( p\)-adic periods, crystalline cohomology and syntomic regulators; \cite{MR1463696}, \cite{MR1145716}, \cite{MR849653} and \cite{MR3264260}. More recently, we are witnessing a growing interest in syntomic and prismatic cohomology, resulting from the work of Antieau, Bhatt, Morrow and Scholze; \cite{MR4502597}, \cite{MR4277271} and \cite{MR4516307}.

Our main result is the following.
\smallskip

\begin{theorem} 
\label{thm:main}
Let \( \Phi\) be a reduced, irreducible root system and   \(G= G_\Phi \) denote the associated simply-connected Chevalley group  over the base field \( \Bbbk\).  If \( \Phi \) is one of the classical systems of type \( \cala_n, \calb_n, \calc_n \) or\, \( \cald_n\), the \lulu\ map \(\mu  \colon \bbl \times \bbu  \times \bbl  \times \bbu  \longrightarrow G  \) is syntomic and faithfully flat.
\end{theorem}

To simplify notation, we henceforth denote \( \bbx := \bbl \times \bbu \times \bbl \times \bbu \). Since the multiplication morphism is always of locally finite presentation, the result above states that that  \( \mu \colon \bbx \to G \) is a faithfully flat, local complete intersection morphism.

In the context of the aforementioned applications (with \( \Bbbk= \bbc\) or \( \bbr\)), this result implies that the Hausdorff measure of the set of  choices for a decomposition of the type \( A = L_1 U_1 L_2 U_2\) is the smallest possible and does not depend on the matrix \( A\).

The genesis of this work stems from applications to algebraic \( K\)-theory. Whenever \( R\) is a regular 
\( \Bbbk\)-algebra, one can use Karoubi-Villamayor \(K\)-theory \( KV_*(R) \) as a model for the algebraic \(K\)-theory of \(R\). Using the results and techniques in \cite{CC-HCG}, one can construct explicit representatives for the Chern classes \( \mbc_{p,n}(\alpha) \in CH^p(\spec{R}, n-2p)\) using homotopy-invariance to identify \(CH^p(\bbx_R,n-2p) \cong  CH^p(\spec{R}, n-2p) \), since \( \bbx \) is \( A^1\)-contractible.  Applying Theorem \ref{thm:main} we can show that any class in \( KV_j(R) \) can be represented by a morphism in \( \sch_{\! /\Bbbk}\)
\[
\alpha \colon \bbx_R \times \Delta^j \longrightarrow SL_{r+1}^{\times j},
\]
for \( r \) sufficiently large, whose restriction to \( \bbx_R \times (\Delta^j - \partial \Delta^j) \) is faithfully flat over \( \Delta^p - \partial \Delta^p \). Such representative \( \alpha \)  is then used to pull-back the   universal cycles \( \mathfrak{C}^p_{j,r+1} \) from \cite{CC-HCG}, yielding explicit representatives of the Chern classes \( \mbc_{p,n}(\alpha)\).

\subsubsection{Notation}
\label{subsubsec:Notation}
Let \( \frg\) denote the Lie algebra of the simply-connected Chevalley group \( G\) over the base field \( \Bbbk\). Fix an order for its roots \( R\) (in the sense of \cite{MR3616493}) and write \( R = R^+ \amalg R^-\), where \( R^+, R^-\) denote the sets of positive and negative roots, respectively. If \( \frh\) is a Cartan subalgebra then
\begin{equation}
\label{eq:LHU}
\frg \ = \ \frl \ \oplus \ \frh \ \oplus \ \fru,
\end{equation}
where \( \fru = \bigoplus_{\alpha \in R^+} \frg_\alpha = \mathsf{Lie}(\bbu) \)\ and\  \( \frl = \bigoplus_{\alpha \in R^+} \frg_{-\alpha}= \mathsf{Lie}(\bbl) \).  As usual, we fix for each root \( \alpha \in R^+\),  elements \( X_\alpha \in \frg_\alpha\), \( Y_\alpha \in \frg_{-\alpha} \) and \( H_\alpha \in \frh\) such that the Lie subalgebra \( \frs_\alpha = \la X_\alpha, Y_\alpha, H_\alpha \ra \)  is isomorphic to \( \frs \frl_2 \).
\medskip

This paper is structured as follows. 
In Section \ref{sec:fibers}, we begin with  an alternative characterization of \( \mu^{-1}(g_0) \), for \( g_0 \in G \), which identifies this fiber with a subscheme \( \scrV_{g_0} \subset \bbl\times \bbu \). We then construct a  morphism \( \varphi \colon \bbu \to  \HomSch{\frl}{\frh} \) and use it to pull-back a determinantal subvariety \( \calb\) of \( \HomSch{\frl}{\frh} \) to \( \bbl\times \calb \subset \bbl \times \bbu \). Later, in Section~\ref {subsec:proof}, we show that an irreducible component of  the fiber \( \mu^{-1}(g_0) \) is contained in the singular set of \(\mu \) only when the corresponding component of  \( \scrV_{g_0} \) is contained in \( \bbl \times \calb.\) We conclude the section by establishing that the subscheme \( \calb\) is defined by a monomial ideal \( \cali_\calb\), as a consequence of the invariance of \( \calb \) under the action of the maximal torus on \( \bbu \) by conjugation.  The main estimate on the dimension of the fibers will follow from the study of the ideal \( \cali_\calb\).

In Section \ref{sec:Cox} we study the nerve complex \( \Delta(\scrA) \) associated to a Coxeter arrangement \( \scrA\).  For simplicity, we denote \(\Delta(\frg):= \Delta(\scrA(\frg) )\)   when \( \scrA(\frg) \) is the arrangement associated to a simple Lie algebra \( \frg\).
Utilizing the well-known structure of the lattice \( L(\scrA(\frg)) \) of the arrangement \(\scrA(\frg)\),  when \( \frg\) is one the classical Lie algebras, we  prove the following result,  where 
\( I_{\Delta(\frg)} \) is the Stanley-Reisner ideal associated to the nerve complex \( \Delta(\frg)\).
\smallskip

\noindent{\bf Theorem \ref{thm:SRdim}.}{\ \it
If \( \frg \) is a classical simple Lie algebra of rank \( n\geq 2\), and \( \Delta(\frg) \) is the nerve complex of its Coxeter arrangement, then the dimension of the corresponding Stanley-Reisner ring \( \Bbbk[\Delta(\frg)] \) is displayed in the following table.
\begin{table}[hpt]
\begin{center}
\begin{tabular}{|c|c|c|}\hline
Type & \( \dim \Bbbk[\Delta(\frg)]\) & \( \operatorname{codim} I_{\Delta(\frg)} \) \\ \hline \hline
\( A_n\) & \( \binom{n}{2}\) & \( n \) \\
\( B_n\) & \( (n-1)^2 \)& \( 2n -1\) \\
\( C_n\) & \( (n-1)^2\) & \( 2n -1\) \\
\( D_n\) & \( 2\binom{n-1}{2} \)   & \( 2(n -1)\)\\ \hline
\end{tabular}
\end{center}
\caption{Dimension and codimension of Stanley-Reisner ideals.}
\end{table}}

In Section \ref{sec:equidim}, we establish a connection between the ideals \( I_{\Delta(\frg)} \) and \(\cali_\calb\).\smallskip

\noindent{\bf Theorem \ref{thm:inclusion}.}{\ \it
Let \( \frg \) be a semi-simple Lie algebra of rank \ \( n \) over \( \Bbbk \). Using the notation above, we have an inclusion of ideals \ 
\( I_{\Delta(\frg)} \ \subseteq \  \cali_\calb \subset \ \Bbbk[\mbx] . \) }
\smallskip

\noindent The proof of  our main result, Theorem \ref{thm:main}, follows from Theorems \ref{thm:SRdim} and \ref{thm:inclusion}.
\smallskip

We conclude with Section \ref{sec:An}, where we present an additional result in the \( \cala_n\) case and pose a few open questions of independent interest.  The main result in this section is the following.
\smallskip

\noindent{\bf Theorem  \ref{prop:mon-An}}{\ \it
For the semi-simple Lie algebra  \( \frs\frl_{n +1} \),  the following equality  holds:
 \[
  I_{\Delta(\frs\frl_{n+1})} \ =  \  \cali_\calb . 
 \]}

In the  \( \cala_n\) case  we also show that the fibers of \lulu\   are global complete intersections, by exhibiting explicit generators for the ideals of each fiber. However, we could not ascertain whether the generators we provide form a regular sequence in any order. The question whether the fibers are global complete intersections in the other cases remain open, and so does the proof of Theorem \ref{prop:mon-An} in the other cases.  

Finally, we observe that the cardinality of the minimal set of monomial generators in the case \( \scrA_n\), for the ideals \( I_{\Delta(\frg)} = \cali_\calb \), is the number of spanning trees of  \( \calk_{n+1} \), the complete graph on \( [n+1] \).   This number is precisely \( \mba(n) := (n+1)^{n-1} \), and this sequence of numbers has many other combinatorial interpretations in the literature. 
In a letter \cite{Del-Loo74} to E. Looijenga, Pierre Deligne calculated\textemdash among other things\textemdash the number of distinct ways one can express a Coxeter element \( c \in W\)  as a product of  \( n = \operatorname{rank}(D) \) reflections, and he denoted this number by \( F'(D)\). It turns out that,  in the \( \scrA_n\) case, this number is precisely \( \mba(n) \).


 


\section{Fibers and critical points of \lulu}
\label{sec:fibers}

This section starts with an alternative description of the fibers of \( \mu \colon \bbx   \to G\), and then proceeds to study its critical points. A key ingredient in this study is a particular map \( \varphi \colon \bbu \to \HomSch{\frl}{\frh} \), see \eqref{eq:varphi}, which is subsequently used to define the closed subscheme \( \calb := \varphi^{-1}\left( D_{\frl, \frh} \right) \subset \bbu \) obtained as the inverse image of a determinantal subvariety \( D_{\frl, \frh} \subset \HomSch{\frl}{\frh}\). The subscheme \( \calb \) plays an essential role  in Section \ref{sec:equidim}.

\subsection{Characterizing the fibers of \( \mu \colon \bbx   \to G\) }

Recall that the multiplication map \( G \times G \to G \) induces closed embeddings 
\begin{equation}
\label{eq:embed}
\begin{tikzcd}
\bbl\times \bbu \arrow[rr, hook, "\phi_-"]  & & G & & \bbu \times \bbl \arrow[ll, hook',"\phi_+"'].
\end{tikzcd}
\end{equation}
\begin{definition}
\label{def:Vg}
Let  \( \bbl\bbu\subset G \)  and  \( \bbu\bbl \subset G \) denote the images of \( \phi_- \) and \( \phi_+\), respectively.
Given  \( g_0\in G\),   define the  subscheme
\begin{equation}
\label{eq:Vg0}
\scrV_{\! g_0} := \ \phi_-^{-1}\left( g_0 \bbu\bbl \right)\   \subset \ \bbl \times \bbu.
\end{equation} 
\end{definition}

\begin{proposition}
\label{prop:iso}
Fix  \( g_0 \in G \), and 
let \( \pi_{12} \colon  \bbx = \bbl\! \times \bbu\!\times \!\bbl\! \times \!\bbu \to \bbl \times \bbu \) denote the projection onto the  first two factors.   Then the restriction of \( \pi_{12} \) to  the fiber \( \mu^{-1}(g_0) \) induces an isomorphism
\(
\pi_{12|_{\mu^{-1}(g_0)}} \ \colon \ \mu^{-1}(g_0) \xrightarrow{ \ \ \cong \ \ } \scrV_{\! g_0}.
\)
\end{proposition}
\begin{proof}
By definition, \( \mu^{-1}(g_0) \) is the \(k(g_0)\)-scheme \( \bbx \times_G \spec{k(g_0)}\). Hence, it suffices to show that,
for any \( k(g_0) \)-algebra \( R\), the projection  \( \pi_{12} \) induces an isomorphism between the corresponding sets of \( R\)-valued points \( \mu^{-1}(g_0)(R) \) and \(\scrV_{\! g_0}(R) \). 

By definition, if  \( (\mba, \mbb) \in \scrV_{\! g_0}(R) \) is given, then \(  \phi_-(\mba,\mbb)=\mba\mbb \) lies in \( (g_0\bbu\bbl) (R) = g_0\cdot \bbu(R)\cdot \bbl(R).\)  Therefore, one can find \( (\mbc, \mbd) \in \bbl(R)\times \bbu(R) \) such that \( \mba\mbb = g_0 \mbd\mbc \), hence
\( \mba\mbb\mbc^{-1}\mbd^{-1} = g_0 \). In other words, \( \mbp:= (\mba, \mbb, \mbc^{-1}, \mbd^{-1}) \in \bbx(R) \) lies in \( \mu^{-1}(g_0)  \) and satisfies \( \pi_{12}(\mbp) = (\mba, \mbb)\). It follows that \( \pi_{12}(R) \colon \mu^{-1}(g_0)(R) \to\scrV_{\! g_0}(R) \) is onto.

On the other hand, suppose that \( \mbp, \mbp' \in \mu^{-1}(g_0)(R) \) satisfy \( \pi_{12}(\mbp)=\pi_{12}(\mbp')\), so that we can write
\( \mbp = (\mba, \mbb, \mbc, \mbd) \) and \( \mbp' =  (\mba, \mbb, \mbc', \mbd') \). Since \( g_0 =\mba\mbb\mbc\mbd=\mba\mbb\mbc'\mbd'\),
one concludes that \( \mbc^{-1}\mbc' = \mbd \mbd'^{-1}  \in \bbu(R) \cap \bbl(R) = \{ 1 \} \). Hence, \( \pi_{12}(R)\) is injective. The result follows. 
\end{proof}

\subsection{The critical points of  $ \mu_{\bar \Bbbk} $ }
\label{subsec:critical}

\newcommand{\geo}[1]{{#1}_{\bar\Bbbk}}

Fix an algebraic closure \( \bar \Bbbk\) of \( \Bbbk\). To simplify notation, given  \( Y \in \sch_{\! / \Bbbk} \), we identify the closed points of  \( \geo{Y}:= Y\times_{\Bbbk} \bar\Bbbk \) with the \(\bar{\Bbbk}\)-valued points \( Y(\bar \Bbbk) \) of \( Y \) over \( \Bbbk\), and 
denote the corresponding base extension of \lulu\ by 
\( \geo{\mu}  \colon \geo{\bbx}  \to \geo{G}  \).

Let \(\geo{\bbx}^\circ \subset \geo{\bbx}  \) be the open subset\textemdash possibly emtpy\textemdash where \( \geo{\mu} \) is smooth, and let  \( \scrC := \geo{\bbx} - \geo{\bbx}^\circ \) denote its complement, with the induced reduced closed subscheme structure, which we call the critical points of \lulu\ (over \(\bar\Bbbk\)).

Through the rest of this section, we   fix a closed point \( p_0 = (x_0, u_0, y_0, v_0) \in \geo{\bbx} \) and denote \( g_0 := \geo{\mu}(p_0) \in \geo{G}\). 
Let 
\( d\mu_{p_0} \colon T_{p_0}\geo{\bbx} \to T_{g_0}\geo{G} \) be the differential of \( \geo{\mu} \) at \( p_0\), and recall that \( \geo{\mu}\) is smooth at \( p_0 \) whenever \( d\mu_{p_0} \) is surjective; see the proof of \cite[Prop. 10.4(iii)]{MR463157}.

Using the identification
\( T_{p_0} \geo{\bbx} = T_{x_0} \geo{\bbl} \oplus T_{u_0} \geo{\bbu} \oplus T_{y_0} \geo{\bbl} \oplus  T_{v_0} \geo{\bbu} \),
 write
\begin{align}
\label{eq:dmu}
d\mu_{p_0}&(A,B,C,D) \\ & = 
d\rho_{u_0y_0v_0}(A) \ + \
d\lambda_{x_0}\circ d\rho_{y_0v_0}(B) \ + \
d\lambda_{x_0u_0}\circ d\rho_{ v_0}(C) \ + \
d\lambda_{x_0u_0y_0}(D), \notag
\end{align}
where  \( \rho_a\) and \( \lambda_a\) denote the morphisms given by right and left translation by \( a \in \geo{G} \), respectively. 
Then, after appropriate right translations identifying tangent spaces with Lie algebras, we obtain a commutative diagram
\begin{equation}\
\label{eq:diff}
\begin{tikzcd}
T_{x_0} \geo{\bbl} \oplus T_{u_0} \geo{\bbu} \oplus T_{y_0} \geo{\bbl} \oplus  T_{v_0} \geo{\bbu}
 \arrow[rr, "d\mu_{p_0}"] \ar[d, "d\rho_{x_0^{-1}}\oplus d\rho_{u_0^{-1}}\oplus d\rho_{y_0^{-1}}\oplus d\rho_{v_0^{-1}}  "', "\cong"] & & T_{g_0}\geo{G}  \arrow[d, "\cong"', "d\rho_{g_0^{-1}}"]\\ 
 \geo{\frl} \oplus  \geo{\fru} \oplus  \geo{\frl} \oplus  \geo{\fru} \arrow[rr, "\widehat{L}_{p_0}"'] & &  \geo{\frg} ,
 \end{tikzcd}
\end{equation}
where \( \widehat{L}_{p_0}\) is defined by
\begin{equation}
\label{eq:hatL}
\widehat{L}_{p_0}(A,B,C,B) \ = \  A + \mathsf{Ad}_{x_0}(B) + \mathsf{Ad}_{x_0u_0}(C)+ \mathsf{Ad}_{x_0u_0y_0}(D).
\end{equation}

\begin{remark}
\label{rem:Hom}
Given finite dimensional \(\Bbbk\) vector spaces \(V \) and \( W\), let
\( \HomSch{V}{W}\) be the affine space 
\( \spec{ \mathsf{Sym_\bullet(\{ V^*\otimes_{\Bbbk} W \}^* } }\) whose \( R \)-valued points  \( \HomSch{V}{W}(R)\)  correspond to \( \Hom{R-\text{mod}}{V_R}{W_R} \), where \( V_R := V\otimes_{\Bbbk} R\) and \( W_R := W\otimes_{\Bbbk} R\). In particular, any morphism of \( \Bbbk\)-schemes \( f \colon X \to \HomSch{V}{W} \), corresponds to a \( \Bbbk\)-linear map  \( f^\sharp \colon V\otimes_\Bbbk W^* \to  \Gamma(X, \scrO_X) \). Equivalently, it corresponds to an element \( f^\flat \in V^*\otimes_\Bbbk W \otimes_\Bbbk \Gamma(X, \scrO_X)\).
\end{remark}

Now, let \( \mathsf{proj}_{\frh} \colon \frg = \frl \oplus \frh \oplus \fru \to \frh  \) denote the linear projection, and define a morphism 
\( \varphi \colon \bbu \to \HomSch{\frl}{\frh}, \ \ u_0 \mapsto \varphi_{u_0}, \) as the composition
\begin{equation}
\label{eq:varphi}
\begin{tikzcd}
\bbu \arrow[r, hook] \arrow[rrrd, bend right=10, "\varphi"']& 
G \arrow[r, "\mathsf{Ad}"] & 
GL(\frg) \arrow[r, "\mathsf{res}_{\frl}"] & 
\HomSch{\frl}{\frg} \arrow[d, "\mathsf{proj}_{\frh *}"] \\
& & & \HomSch{\frl}{\frh}.
\end{tikzcd}
\end{equation}

\begin{definition}
\label{def:varphi2}
Given \( p_0 = (x_0, u_0, y_0, v_0)\in  \bbx(\bar{\Bbbk}) \)\  
define 
\begin{align}
\label{eq:Lp0}
 L_{p_0} \colon \geo{\frl}\oplus \geo{\fru}\oplus \geo{\frl}\oplus \geo{\fru} &\ \  \longrightarrow \ \  \geo{\frg}  \\
(A,B,C,D)& \ \ \longmapsto \ \ A + B + \Ad{u_0}{C} + \Ad{u_0y_0}{D}. \notag
\end{align}
\end{definition}

\begin{proposition}
\label{prop:varphi3}
With  \( p_0  = (x_0, u_0, y_0, v_0)\in  \bbx(\bar{\Bbbk}) \)\  and\  \( g_0 = \geo{\mu}(p_0) \in  {G}(\bar{\Bbbk}) \) as above, the following holds. 
\begin{enumerate}[i.]
\item The differential \( d\mu_{p_0}\colon T_{p_0}\geo{\bbx} \to T_{g_0}\geo{G}  \) at \( p_0 \) is onto if and only if \( L_{p_0} \) is onto. 
\item If \( \varphi_{u_0} \colon \geo{\frl} \to \geo{\frh} \) is onto, then so is \( d\mu_{p_0} \).
\end{enumerate}
\end{proposition}
\begin{proof}
First observe that \( \mathsf{Ad}_{x_0} \) restricts to an automorphism of \( \geo{\frl}\), since \( x_0 \in \bbl_{\bar\Bbbk}\). Therefore, one can write \( L_{p_0}\) as the composition
\begin{equation}
\begin{tikzcd}
\geo{\frl}\oplus \geo{\fru}\oplus \geo{\frl}\oplus \geo{\fru}
\arrow[rrrrrd, "L_{p_0}", bend right=10]
 \arrow[rrr, "\mathsf{Ad}_{x_0}\oplus I \oplus I \oplus I"', "\cong"] & & &  \geo{\frl}\oplus \geo{\fru}\oplus \geo{\frl}\oplus \geo{\fru} \ar[rr, "\widehat{L}_{p_0}"'] & & \geo{\frg} \arrow[d, "\cong"', "\mathsf{Ad}_{x^{-1}_0}"] \\
& & & & &  \geo{\frg},
\end{tikzcd}
\end{equation}
where \( \widehat{L}_{p_0} \) is defined in \eqref{eq:hatL}.  Hence, \( L_{p_0} \) is onto if and only if so is \( \widehat L_{p_0} \), and the proof of the first assertion follows from the commutativity of the diagram in \eqref{eq:diff}.

Now, consider \( V = A + H + B \ \in\ \geo{\frg}  =  \geo{\frl} \oplus \geo{\frh} \oplus \geo{\fru} \). If \( \varphi_{u_0} \) is onto,  one can find \( C \in \geo{\frl} \) such that \( \Ad{u_0}{C} = A' + H + B' \), for some \( A' \in \geo{\frl} \) and \( B' \in \geo{\fru}\).  Then \( L_{p_0}(A-A', B-B', C, 0) = A-A'+ B-B' + A' + H + B' = A+H+B = V\) . In other words,  \( L_{p_0} \) is surjective, and the second assertion in the proposition now follows from the first one. 
\end{proof}

\subsection{Pulling back determinantal varieties}
\label{subsec:pb}

Let \( D_{\frl, \frh} \subset \HomSch{\frl}{\frh}\) be the determinantal variety over \( \Bbbk\) whose set of \( R \)-valued points \( D_{\frl, \frh}(R) \), for    \(\Bbbk\subset  R\), consists of those \( R\)-linear maps from \( \frl_R\) to \( \frh_R\) which do not have maximal rank. Define
\begin{equation}
\label{eq:B}
\calb \ := \ \varphi^{-1}(D_{\frl,\frh}) \ \subset \ \bbu. 
\end{equation}
Here we show that the subscheme \(\calb\) is defined by a monomial ideal, and in the next section we provide estimates for its codimension in \( \bbu\). 

\begin{proposition}
\label{prop:torus-inv}
The subscheme \( \calb\) is invariant under the conjugation action of the torus on \( \bbu\).
\end{proposition}
\begin{proof}
Consider the composition

\[
\begin{tikzcd}
\bbt \times_\Bbbk \calb \arrow[r, hook] \ar[rrr, rounded corners=8pt, to path={-- ([yshift=-2.5ex]\tikztostart.south)-| (\tikztotarget)}]  &\bbt \times_\Bbbk \bbu \arrow[r, "\widetilde{\mathsf Ad}"] \ar[d, phantom, "\scriptstyle \widetilde\varphi", pos=0.85]  &  \bbu \arrow[r, "\varphi"] & \HomSch{\frl}{\fru},\\
& \phantom{.} & & 
\end{tikzcd}
\]
where \( \widetilde{\mathsf Ad} \) is the adjoint of \( \mathsf{Ad} \colon \bbt \to \mathsf{Aut}(\bbu) \). Given a \(\Bbbk\)-algebra \( R \) and  \( ( t,u_0) \in \bbt(R)\times \bbu(R) \), we just need to show that the  morphisms \(\widetilde{\varphi}( t,u_0)  \colon \frl_R \to \fru_R\)
and   \(\widetilde{\varphi}( 1,u_0)  \colon \frl_R \to \fru_R\) have the same image, where \( 1 \in \bbt(R)\) is the identity element.

First note that for \( C \in \frl_R \) and \( t \in \bbt(R)\), the element
\( C' := \Ad{t^{-1}}{C} \) also lies in \( \frl_R\), since \(\mathsf{Ad}_{t^{-1}} \) preserves the decomposition \( \frg = \frl\oplus \frh \oplus \fru \), for all \(t \in \bbt \). Now, write
\(
\Ad{u_0}{C'} \ = \ A \oplus \varphi_{u_0}(C') \oplus B,
\)
with \( A \in \frl_R\) and \(B \in \fru_R\). Then
\begin{equation}
\label{eq:u0}
\Ad{t}{\Ad{u_0}{C'}} \ = \ \Ad{t}{A} + \varphi_{u_0}(C') + \Ad{t}{B},
\end{equation}
and hence
\begin{align*}
\widetilde{\varphi}(t, u_0)(C) 
&= 
\mathsf{proj}_\frh\left( \Ad{tu_0t^{-1}}{C} \right)    = 
\mathsf{proj}_\frh\left( \mathsf{Ad}_{t}\circ \mathsf{Ad}_{u_0}\circ \mathsf{Ad}_{t^{-1}}(C) \right)  \\
& =
\mathsf{proj}_\frh\left( \mathsf{Ad}_{t}\circ \mathsf{Ad}_{u_0}(C') \right)  =
\varphi_{u_0}(C'),
\end{align*}
where the last identity follows from \eqref{eq:u0}. This suffices to show the result.
\end{proof}

For each root \( \alpha \in R \), let \( \lambda_\alpha \colon \bbg_a \to G \) denote the associated unipotent element. Fix a  listing \(R^+= \{ \alpha_1, \alpha_2, \ldots, \alpha_D \} \) of the positive roots to define an isomorphism
\begin{align}
\label{eq:iso}
\Psi \colon \bba^D & \xrightarrow{\ \cong\ } \bbu \\
(x_1, x_2, \ldots, x_D) & \longmapsto \lambda_{\alpha_1}(x_1)\, \lambda_{\alpha_2}(x_2) \cdots \lambda_{\alpha_D}(x_D). \notag
\end{align}
We often use the symbol \( x_\alpha\) to denote the variable \( x_i \) corresponding to the label  of \( \alpha \in R^+\) in this listing. Hence,  
by construction, the differential of \( \Psi \) at \( 0 \in \bba^D\) 
\begin{equation}
\label{eq:diffPsi}
d\Psi_0 \colon T_0\bba^D \longrightarrow T_I\bbu = \fru 
\end{equation}
satisfies \( d\Psi_0 (\frac{\partial}{\partial x_\alpha})  = X_\alpha, \)  for each \( \alpha\in R^+\); see Notation \ref{subsubsec:Notation}.

Now, let \( \Bbbk[\mbx] \) denote the polynomial ring \( \Bbbk[x_1, \ldots, x_D] \)  and denote by \( \cali_\calb \subset  \Bbbk[\mbx]\) the ideal of the closed subscheme \( \Psi^{-1}(\calb) \subset \bba^D\).

\begin{corollary}
\label{cor:mon}
\( \cali_\calb\) is a monomial ideal. 
\end{corollary}
\begin{proof}
Suppose \( p(\mbx) = \sum_M\, c_M \, \mbx^M \ \in \cali_\calb \), where \( M = (m_1, \ldots, m_D) \in \bbn^D \) and 
\( \mbx^M:= x_1^{m_1} \cdots x_D^{m_D} \).  Therefore, if \( \mba^0 \in \Psi^{-1}(\calb)(K) \) and \( t \in \bbt(K)\) with \( \Bbbk\subset K\), then
\( t \Psi(\mba^0)t^{-1} \) lies in \( \calb(K)\), according to the proposition. 
However, 
\begin{align*}
 t \Psi(\mba^0)t^{-1} & = 
 t \left\{ \prod_{\ell=1}^D\, \lambda_{\alpha_\ell}(a^0_\ell)     \right\}t^{-1} \ = \
  \prod_{\ell=1}^D\, \left\{ t \, \lambda_{\alpha_\ell}(a^0_\ell)\,     t^{-1} \right\} \\
  & =
  \prod_{\ell=1}^D\,     \lambda_{\alpha_\ell}(\alpha_\ell(t)\, a^0_\ell) \  =: \ \Psi(t\star \mba^0) ,
\end{align*}
where \( t \star \mba := (\alpha_1(t)\, a^0_1, \ldots, \alpha_D(t)\, a^0_D) \). 

Therefore, if \( p(\mbx) \) lies in \( \cali_\calb\) so does \( p(t\star \mbx) \), for all \( t\in \bbt(\Bbbk) \). 
However \( p(t\star~\mbx) = \sum_M\, c_M \left( \prod_{\ell=1}^D \alpha(t)^{m_\ell} \right)\mbx^M\), and from the linear independence of characters one concludes that \( c_M \mbx^M \) lies in \(\cali_\calb\) for each \( M \) such that \( c_M \neq 0\), and so does \( \mbx^M\).
\end{proof}

\subsubsection{Generators for \(\cali_\calb\)}
\label{subsec:gens}

Fix a basis \( \{ \mbv_1,\ldots, \mbv_n\} \) for \( \frh \), and let \( \{ Y^*_{\alpha_1}, \ldots, Y^*_{\alpha_D} \} \subset \frl^*  \) be the basis dual to \( \{ Y_{\alpha_1}, \ldots, Y_{\alpha_D} \} \). 
Using Remark \ref{rem:Hom}, we identify the composition
\[
\begin{tikzcd}
\bba^D_{  \Bbbk}  \arrow[rr, "\Psi", "\cong"'] & & \bbu \arrow[rr, "\varphi"] & & \HomSch{\frl }{\frh }
\end{tikzcd}
\] 
with an element \( \mbp(\mbx) \in \frl^* \otimes \frh \otimes \Bbbk[\mbx] \) that can be uniquely written as 
\[
\mbp(\mbx) \ = \ \sum_{i=1}^D \sum_{j=1}^n Y^*_{\alpha_i}\otimes \mbv_j \otimes p_{ij}(\mbx) .
\]

Given \( I = (i_1<\cdots < i_n) \in \Lwedge^n[D] \), define
\begin{equation}
\label{eq:qI}
q_I(\mbx) \ := \ \det \left(   p_{i_r, s}(\mbx)  \right)_{1\leq r, s \leq n }  \ \in \ \bar \Bbbk[\mbx].
\end{equation}
Then, by definition, \(\cali_\calb = \la q_I(\mbx) \mid I \in \Lwedge^n[D] \ra . \)

By construction, \( \varphi  \circ \Psi (0) = 0 \), and one concludes that \( q_I(\mbx) \in \frakm^n \), for all \( I \in \Lwedge^n[D] \), where \( \frakm \) is the maximal ideal 
\( \la x_1, \ldots, x_D\ra \subset   \Bbbk[\mbx] \).
Now, write \( q_I(\mbx) = \sum_{M\in \calk_I}\, c_M^I \mbx^M\), where
\( \calk_I := \{ M\in \bbn^D \mid c^I_M \neq 0 \} = \operatorname{supp}(q_I(\mbx)) \) is the support of the polynomial \( q_I(\mbx) \). 
It follows from Corollary \ref{cor:mon} that one has a presentation
\begin{equation}
\label{eq:caliB}
\cali_\calb = \biggl< \mbx^M \ \mid \  M\  \in \bigcup_{I \in \Lwedge^n[D]}\, \calk_I  \biggr> .
\end{equation}



\section{Coxeter arrangements and monomial ideals}
\label{sec:Cox}

In this section we make a little digression to study the nerve complex  \( \Delta(\scrA)\) associated to a Coxeter arrangement \( \scrA\). The main objective is to determine the codimension of the corresponding Stanley-Reisner ideal \( I_{\Delta(\scrA)} \) in each of the classical cases. Although there is no claim to originality in this section, we were unable to find direct references in the literature to the statements in Theorem \ref{thm:SRdim}.

\subsection{The nerve complex}
\label{subsec:nerve}

Let \( V \) be a real vector space, and consider  a central hyperplane arrangement \( \scrA = \{ \Omega_1, \ldots, \Omega_D\} \subset V .\)  If \( \left( L(\scrA), \preceq \right)\) is the intersection lattice of \( \scrA\), let   \( \widehat 0  := V\) and \( \widehat \bone := \cap_{\ell=1}^D \Omega_\ell \) denote its minimal and maximal elements, respectively.
Denote \( L(\scrA)^\circ := L(\scrA) - \{ \widehat 0, \widehat \bone \} \) and, for a given \( \pi \in L(\scrA)^\circ \), define 
\begin{equation}
\label{eq:flat}
F_\pi \ := \ \{ \Omega \in \scrA \mid \Omega \prec \pi \}.
\end{equation}

\begin{definition}
\label{def:nerve}
Given \( \scrA\) as above, let \( \Delta(\scrA)\) denote its associated \emph{nerve complex}. This is a simplicial complex having \( \scrA \) as its vertex set, and whose faces are generated by \( \{ F_\pi \mid \pi \in L(\scrA)^\circ \} \).
Let \( I_{\Delta(\scrA)} \subset \Bbbk[x_1, \ldots, x_D] \) denote the corresponding Stanley-Reisner ideal, generated by the non-face monomials, i.e., 
\begin{equation}
\label{eq:SRI}
I_{\Delta(\scrA)} \ = \ \la x_{r_1}\cdots x_{r_k} \ \mid \ \{ \Omega_{r_1}, \ldots, \Omega_{r_k}\} \ \text{ is not a face in } \Delta(\scrA) \ra.
\end{equation}
The Stanley-Reisner ring \( \Bbbk[\Delta(\scrA)] \) is the quotient
\( 
 \Bbbk[x_1, \ldots, x_D]/I_{\Delta(\scrA)}
\); see  \cite{MR1453579}.

\end{definition}

\begin{remark} 
\label{rem:SR}
\begin{enumerate}[a.]
\item The complex \( \Delta(\scrA)\) depends only on the lattice \( L(\scrA)\), and isomorphic lattices yield isomorphic Stanley-Reisner rings.
\item It follows from \cite[II, Thm. 1.3]{MR1453579} that 
\begin{equation}
\label{eq:dimSRR}
\dim  \left(  \Bbbk[\Delta(\scrA)]\right) = \dim \Delta(\scrA) + 1 = 
\max\, \{ \# F_\pi \ \mid\ \pi \in L(\scrA)^\circ \}.
\end{equation}
\end{enumerate}
\end{remark}

The examples of interest in this paper are the \emph{Coxeter arrangements} of type \( A_n, B_n, C_n\) and \( D_n\),  defined as follows. Consider a classical simple Lie algebra (over \( \Bbbk\)) of one of these types, and recall Notation \ref{subsubsec:Notation}.  Then, define the arrangement
\begin{equation}
\label{eq:CoxArrang}
\scrA(\frg) := \{ \Omega_\alpha  \mid  \alpha \in R^+\},
\end{equation}
where, for \( \alpha \in R^+ \subset \frh^* \) we denote
\(
\Omega_\alpha := \{ H \in \frh \mid \la H, \alpha \ra = 0 \}.
\)
In particular, the Weyl group \( W \) is generated by the reflections across the radicial hyperplanes \( \Omega_\alpha \in \scrA(\frg) \).   
In this case, we denote the  complex \( \Delta(\scrA(\frg)) \) simply by \( \Delta(\frg)\).

\begin{example}
\label{exmp:A3}
Consider the \( A_3 \) arrangement, so that \( \frg = \frs\frl_4\). In this case, one identifies 
\( \frh_\bbr \equiv \{ (x_1,x_2,x_3,x_4) \mid x_1 +x_2 +x_3 +x_4 = 0 \} \subset \bbr^4 \equiv \bbr\{ \mbe_1, \mbe_2, \mbe_3, \mbe_4\}. 
\)
Denote by \( \mbe_{ij}^* \in \frh_\bbr^*\) the restriction of \( \mbe_i^* - \mbe_j^* \in \bbr^{4*} \) to \( \frh_\bbr\), and identify the set of positive roots with \( R^+ = \{ \mbe_{ij}^* \mid 1\leq i<j \leq 4 \} \).  The corresponding arrangement is
\[
\scrA = \{ \Omega_{12},  \Omega_{13},  \Omega_{14},  \Omega_{23},  \Omega_{24},  \Omega_{34} \},
\]
where \( \Omega_{ij} = \{ (x_1, x_2, x_3, x_4 ) \mid x_i = x_j \} \).
The intersection lattice \( L(\scrA(\frs\frl_4)) \) is depicted below, while the notation for the ``flats'' in the lattice will be explained in the next section. 

%
\begin{center}
\begin{figure}[h]

\tikzset{every picture/.style={line width=0.75pt}} 

\begin{tikzpicture}[x=0.32pt,y=0.3pt,yscale=-1,xscale=1]

\draw   (200,517.78) -- (198.62,522.02) -- (195.01,524.65) -- (190.55,524.65) -- (186.94,522.02) -- (185.56,517.78) -- (186.94,513.54) -- (190.55,510.92) -- (195.01,510.92) -- (198.62,513.54) -- cycle ;
\draw   (300,517.78) -- (298.62,522.02) -- (295.01,524.65) -- (290.55,524.65) -- (286.94,522.02) -- (285.56,517.78) -- (286.94,513.54) -- (290.55,510.92) -- (295.01,510.92) -- (298.62,513.54) -- cycle ;
\draw   (400,517.78) -- (398.62,522.02) -- (395.01,524.65) -- (390.55,524.65) -- (386.94,522.02) -- (385.56,517.78) -- (386.94,513.54) -- (390.55,510.92) -- (395.01,510.92) -- (398.62,513.54) -- cycle ;
\draw   (500,517.78) -- (498.62,522.02) -- (495.01,524.65) -- (490.55,524.65) -- (486.94,522.02) -- (485.56,517.78) -- (486.94,513.54) -- (490.55,510.92) -- (495.01,510.92) -- (498.62,513.54) -- cycle ;
\draw   (700,517.78) -- (698.62,522.02) -- (695.01,524.65) -- (690.55,524.65) -- (686.94,522.02) -- (685.56,517.78) -- (686.94,513.54) -- (690.55,510.92) -- (695.01,510.92) -- (698.62,513.54) -- cycle ;
\draw   (600,517.78) -- (598.62,522.02) -- (595.01,524.65) -- (590.55,524.65) -- (586.94,522.02) -- (585.56,517.78) -- (586.94,513.54) -- (590.55,510.92) -- (595.01,510.92) -- (598.62,513.54) -- cycle ;
\draw   (100,343.28) -- (96.64,349.1) -- (89.92,349.1) -- (86.56,343.28) -- (89.92,337.46) -- (96.64,337.46) -- cycle ;
\draw   (802.44,343.28) -- (799.08,349.1) -- (792.36,349.1) -- (789,343.28) -- (792.36,337.46) -- (799.08,337.46) -- cycle ;
\draw   (325,343.28) -- (321.64,349.1) -- (314.92,349.1) -- (311.56,343.28) -- (314.92,337.46) -- (321.64,337.46) -- cycle ;
\draw   (575,343.28) -- (571.64,349.1) -- (564.92,349.1) -- (561.56,343.28) -- (564.92,337.46) -- (571.64,337.46) -- cycle ;
\draw   (198.43,349.83) -- (185.67,335.4) -- (210.8,335.06) -- cycle ;
\draw   (452.43,349.83) -- (439.63,332.95) -- (464.77,332.61) -- cycle ;
\draw   (687.43,350.17) -- (674.67,335.73) -- (699.8,335.4) -- cycle ;
\draw [color={rgb, 255:red, 208; green, 2; blue, 27 }  ,draw opacity=1 ]   (93.28,343.28) -- (192.78,517.78) ;
\draw [color={rgb, 255:red, 208; green, 2; blue, 27 }  ,draw opacity=1 ]   (93.28,343.28) -- (292.78,517.78) ;
\draw [color={rgb, 255:red, 208; green, 2; blue, 27 }  ,draw opacity=1 ]   (93.28,343.28) -- (392.78,517.78) ;
\draw [color={rgb, 255:red, 65; green, 117; blue, 5 }  ,draw opacity=1 ]   (318.28,343.28) -- (192.78,517.78) ;
\draw [color={rgb, 255:red, 65; green, 117; blue, 5 }  ,draw opacity=1 ]   (318.28,343.28) -- (492.78,517.78) ;
\draw [color={rgb, 255:red, 65; green, 117; blue, 5 }  ,draw opacity=1 ]   (318.28,343.28) -- (592.78,517.78) ;
\draw [color={rgb, 255:red, 144; green, 19; blue, 254 }  ,draw opacity=1 ]   (568.28,343.28) -- (292.78,517.78) ;
\draw [color={rgb, 255:red, 144; green, 19; blue, 254 }  ,draw opacity=1 ]   (568.28,343.28) -- (492.78,517.78) ;
\draw [color={rgb, 255:red, 144; green, 19; blue, 254 }  ,draw opacity=1 ]   (568.28,343.28) -- (692.78,517.78) ;
\draw [color={rgb, 255:red, 139; green, 87; blue, 42 }  ,draw opacity=1 ]   (795.72,343.28) -- (392.78,517.78) ;
\draw [color={rgb, 255:red, 139; green, 87; blue, 42 }  ,draw opacity=1 ]   (795.72,343.28) -- (592.78,517.78) ;
\draw [color={rgb, 255:red, 139; green, 87; blue, 42 }  ,draw opacity=1 ]   (795.72,343.28) -- (692.78,517.78) ;
\draw [color={rgb, 255:red, 27; green, 76; blue, 131 }  ,draw opacity=1 ] [dash pattern={on 4.5pt off 4.5pt}]  (198.33,342.53) -- (192.78,517.78) ;
\draw [color={rgb, 255:red, 31; green, 90; blue, 158 }  ,draw opacity=1 ] [dash pattern={on 4.5pt off 4.5pt}]  (198.33,342.53) -- (692.78,517.78) ;
\draw [color={rgb, 255:red, 208; green, 2; blue, 27 }  ,draw opacity=1 ] [dash pattern={on 4.5pt off 4.5pt}]  (452.32,341.31) -- (492.78,517.78) ;
\draw [color={rgb, 255:red, 208; green, 2; blue, 27 }  ,draw opacity=1 ] [dash pattern={on 4.5pt off 4.5pt}]  (452.32,341.31) -- (392.78,517.78) ;
\draw [color={rgb, 255:red, 144; green, 19; blue, 254 }  ,draw opacity=1 ] [dash pattern={on 4.5pt off 4.5pt}]  (687.33,342.87) -- (292.78,517.78) ;
\draw [color={rgb, 255:red, 144; green, 19; blue, 254 }  ,draw opacity=1 ] [dash pattern={on 4.5pt off 4.5pt}]  (687.33,342.87) -- (592.78,517.78) ;
\draw   (462.5,125) -- (475,137.5) -- (462.5,150) -- (450,137.5) -- cycle ;
\draw   (437.5,667) -- (450,679.5) -- (437.5,692) -- (425,679.5) -- cycle ;
\draw  [dash pattern={on 0.84pt off 2.51pt}]  (462.5,137.5) -- (93.28,343.28) ;
\draw  [dash pattern={on 0.84pt off 2.51pt}]  (462.5,137.5) -- (198.33,342.53) ;
\draw  [dash pattern={on 0.84pt off 2.51pt}]  (462.5,137.5) -- (318.28,343.28) ;
\draw  [dash pattern={on 0.84pt off 2.51pt}]  (462.5,137.5) -- (452.32,341.31) ;
\draw  [dash pattern={on 0.84pt off 2.51pt}]  (462.5,137.5) -- (568.28,343.28) ;
\draw  [dash pattern={on 0.84pt off 2.51pt}]  (462.5,137.5) -- (687.33,342.87) ;
\draw  [dash pattern={on 0.84pt off 2.51pt}]  (462.5,137.5) -- (795.72,343.28) ;
\draw  [dash pattern={on 0.84pt off 2.51pt}]  (192.78,517.78) -- (437.5,679.5) ;
\draw  [dash pattern={on 0.84pt off 2.51pt}]  (292.78,517.78) -- (437.5,679.5) ;
\draw  [dash pattern={on 0.84pt off 2.51pt}]  (392.78,517.78) -- (437.5,679.5) ;
\draw  [dash pattern={on 0.84pt off 2.51pt}]  (492.78,517.78) -- (437.5,679.5) ;
\draw  [dash pattern={on 0.84pt off 2.51pt}]  (592.78,517.78) -- (437.5,679.5) ;
\draw  [dash pattern={on 0.84pt off 2.51pt}]  (692.78,517.78) -- (437.5,679.5) ;
%
\draw (76,306) node [anchor=north west][inner sep=0.75pt]   [align=left] {\footnotesize 123|4};
\draw (672,306) node [anchor=north west][inner sep=0.75pt]   [align=left] {\footnotesize 13|24};
\draw (435,302) node [anchor=north west][inner sep=0.75pt]   [align=left] {\footnotesize 14|23};
\draw (787,306) node [anchor=north west][inner sep=0.75pt]   [align=left] {\footnotesize 1|234};
\draw (160,531) node [anchor=north west][inner sep=0.75pt]   [align=left] {\footnotesize 12|3|4};
\draw (256,531) node [anchor=north west][inner sep=0.75pt]   [align=left] {\footnotesize 13|2|4};
\draw (356,531) node [anchor=north west][inner sep=0.75pt]   [align=left] {\footnotesize 23|1|4};
\draw (460,531) node [anchor=north west][inner sep=0.75pt]   [align=left] {\footnotesize 14|2|3};
\draw (566,531) node [anchor=north west][inner sep=0.75pt]   [align=left] {\footnotesize 24|1|3};
\draw (666,531) node [anchor=north west][inner sep=0.75pt]   [align=left] {\footnotesize 34|1|2};
\draw (176,306) node [anchor=north west][inner sep=0.75pt]   [align=left] {\footnotesize 12|34};
\draw (301,306) node [anchor=north west][inner sep=0.75pt]   [align=left] {\footnotesize 124|3};
\draw (551,302) node [anchor=north west][inner sep=0.75pt]   [align=left] {\footnotesize 134|2};
\draw (428,698) node [anchor=north west][inner sep=0.75pt]   [align=left] {\(\bf \widehat 0\)};
\draw (452,84) node [anchor=north west][inner sep=0.75pt]   [align=left] {\( \bf \widehat 1\)};
\end{tikzpicture}
\caption{ Intersection lattice of the \(A_3\) arrangement.}
\end{figure}
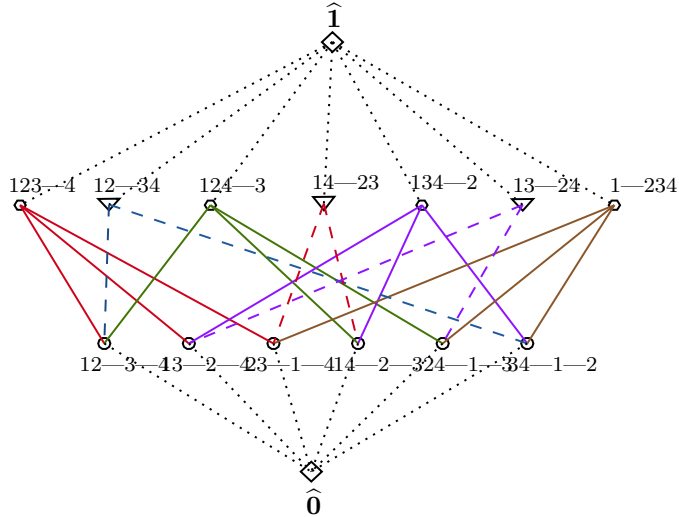
\end{center}

Using the intersection lattice above, one can directly see that the complex  \( \Delta(\scrA_3)\) has four \( 2\)-dimensional facets
\( \left\{ F_{123|4}, F_{124|3}, F_{134|2}, F_{1|234} \right\} \) and three 
\( 1 \)-dimensional facets \( \left\{ F_{12|34}, F_{14|23}, F_{13|24} \right\}\) that assemble as displayed in the next figure.


%
\begin{center}
\begin{figure}[h]

  
\tikzset {_fe20tsehh/.code = {\pgfsetadditionalshadetransform{ \pgftransformshift{\pgfpoint{0 bp } { 0 bp }  }  \pgftransformrotate{0 }  \pgftransformscale{2 }  }}}
\pgfdeclarehorizontalshading{_csuf7xqb2}{150bp}{rgb(0bp)=(1,1,1);
rgb(37.5bp)=(1,1,1);
rgb(62.5bp)=(0.9,0.9,0.9);
rgb(100bp)=(0.9,0.9,0.9)}

  
\tikzset {_bgoua249b/.code = {\pgfsetadditionalshadetransform{ \pgftransformshift{\pgfpoint{0 bp } { 0 bp }  }  \pgftransformrotate{0 }  \pgftransformscale{2 }  }}}
\pgfdeclarehorizontalshading{_lwdhe5o8o}{150bp}{rgb(0bp)=(0,0.5,0.5);
rgb(39.285714285714285bp)=(0,0.5,0.5);
rgb(59.82142857142857bp)=(1,1,0);
rgb(100bp)=(1,1,0)}

  
\tikzset {_flungya80/.code = {\pgfsetadditionalshadetransform{ \pgftransformshift{\pgfpoint{0 bp } { 0 bp }  }  \pgftransformrotate{0 }  \pgftransformscale{2 }  }}}
\pgfdeclarehorizontalshading{_1gw1wh7zt}{150bp}{rgb(0bp)=(1,1,0);
rgb(37.5bp)=(1,1,0);
rgb(62.5bp)=(0,0.5,0.5);
rgb(100bp)=(0,0.5,0.5)}

  
\tikzset {_ff01gojv6/.code = {\pgfsetadditionalshadetransform{ \pgftransformshift{\pgfpoint{0 bp } { 0 bp }  }  \pgftransformrotate{0 }  \pgftransformscale{2 }  }}}
\pgfdeclarehorizontalshading{_otc30b6nt}{150bp}{rgb(0bp)=(1,1,0);
rgb(37.5bp)=(1,1,0);
rgb(62.5bp)=(0,0.5,0.5);
rgb(100bp)=(0,0.5,0.5)}
\tikzset{every picture/.style={line width=0.75pt}} 

\begin{tikzpicture}[x=0.75pt,y=0.75pt,yscale=-1,xscale=1]

\path  [shading=_csuf7xqb2,_fe20tsehh] (260.33,78) -- (368.03,259.04) -- (152.08,259.04) -- cycle ; 
 \draw  [color={rgb, 255:red, 139; green, 87; blue, 42 }  ,draw opacity=1 ][line width=0.75]  (260.33,78) -- (368.03,259.04) -- (152.08,259.04) -- cycle ; 

\draw    (259.41,78) .. controls (261.09,79.65) and (261.11,81.31) .. (259.46,83) .. controls (257.81,84.68) and (257.83,86.35) .. (259.51,88) .. controls (261.19,89.65) and (261.21,91.32) .. (259.56,93) .. controls (257.91,94.68) and (257.93,96.35) .. (259.61,98) .. controls (261.29,99.65) and (261.31,101.32) .. (259.66,103) .. controls (258.01,104.68) and (258.03,106.35) .. (259.71,108) .. controls (261.39,109.65) and (261.41,111.32) .. (259.76,113) .. controls (258.11,114.68) and (258.13,116.35) .. (259.81,118) .. controls (261.49,119.65) and (261.51,121.32) .. (259.86,123) .. controls (258.21,124.68) and (258.23,126.35) .. (259.91,128) .. controls (261.59,129.65) and (261.61,131.32) .. (259.96,133) .. controls (258.31,134.68) and (258.33,136.35) .. (260.01,138) .. controls (261.69,139.65) and (261.71,141.32) .. (260.06,143) .. controls (258.41,144.68) and (258.43,146.35) .. (260.11,148) .. controls (261.79,149.65) and (261.81,151.32) .. (260.16,153) .. controls (258.51,154.68) and (258.53,156.35) .. (260.21,158) .. controls (261.89,159.65) and (261.91,161.32) .. (260.26,163) .. controls (258.61,164.68) and (258.63,166.35) .. (260.31,168) .. controls (261.99,169.65) and (262.01,171.32) .. (260.36,173) .. controls (258.71,174.68) and (258.73,176.34) .. (260.41,177.99) .. controls (262.09,179.64) and (262.11,181.31) .. (260.46,182.99) .. controls (258.81,184.68) and (258.83,186.34) .. (260.52,187.99) .. controls (262.2,189.64) and (262.22,191.31) .. (260.57,192.99) .. controls (258.92,194.67) and (258.94,196.34) .. (260.62,197.99) .. controls (262.3,199.64) and (262.32,201.31) .. (260.67,202.99) .. controls (259.02,204.67) and (259.04,206.34) .. (260.72,207.99) .. controls (262.4,209.64) and (262.42,211.31) .. (260.77,212.99) .. controls (259.12,214.67) and (259.14,216.34) .. (260.82,217.99) .. controls (262.5,219.64) and (262.52,221.31) .. (260.87,222.99) .. controls (259.22,224.67) and (259.24,226.34) .. (260.92,227.99) .. controls (262.6,229.64) and (262.62,231.31) .. (260.97,232.99) .. controls (259.32,234.67) and (259.34,236.34) .. (261.02,237.99) .. controls (262.7,239.64) and (262.72,241.31) .. (261.07,242.99) .. controls (259.42,244.67) and (259.44,246.34) .. (261.12,247.99) .. controls (262.8,249.64) and (262.82,251.31) .. (261.17,252.99) .. controls (259.52,254.67) and (259.54,256.34) .. (261.22,257.99) .. controls (262.9,259.64) and (262.92,261.31) .. (261.27,262.99) .. controls (259.62,264.67) and (259.64,266.34) .. (261.32,267.99) .. controls (263,269.64) and (263.02,271.31) .. (261.37,272.99) .. controls (259.72,274.67) and (259.74,276.34) .. (261.42,277.99) .. controls (263.1,279.64) and (263.12,281.31) .. (261.47,282.99) .. controls (259.82,284.67) and (259.84,286.34) .. (261.52,287.99) .. controls (263.2,289.64) and (263.22,291.31) .. (261.57,292.99) .. controls (259.92,294.67) and (259.94,296.34) .. (261.62,297.99) .. controls (263.3,299.64) and (263.32,301.31) .. (261.67,302.99) .. controls (260.02,304.68) and (260.04,306.34) .. (261.73,307.99) -- (261.76,311.67) -- (261.76,311.67) ;
\draw  [dash pattern={on 4.5pt off 4.5pt}]  (152.08,259.04) -- (369,132.67) ;
\draw  [dash pattern={on 0.75pt off 0.75pt}]  (149.81,133.13) .. controls (152.08,132.52) and (153.53,133.36) .. (154.14,135.63) .. controls (154.75,137.91) and (156.19,138.74) .. (158.47,138.13) .. controls (160.75,137.52) and (162.19,138.35) .. (162.8,140.63) .. controls (163.41,142.91) and (164.86,143.74) .. (167.14,143.13) .. controls (169.42,142.52) and (170.86,143.35) .. (171.47,145.63) .. controls (172.08,147.91) and (173.52,148.74) .. (175.8,148.13) .. controls (178.08,147.52) and (179.52,148.35) .. (180.13,150.63) .. controls (180.74,152.9) and (182.19,153.73) .. (184.46,153.12) .. controls (186.74,152.51) and (188.18,153.34) .. (188.79,155.62) .. controls (189.4,157.9) and (190.84,158.73) .. (193.12,158.12) .. controls (195.4,157.51) and (196.84,158.34) .. (197.45,160.62) .. controls (198.06,162.9) and (199.5,163.73) .. (201.78,163.12) .. controls (204.06,162.51) and (205.5,163.34) .. (206.11,165.62) .. controls (206.72,167.9) and (208.16,168.73) .. (210.44,168.12) .. controls (212.72,167.51) and (214.16,168.34) .. (214.77,170.62) .. controls (215.39,172.89) and (216.84,173.72) .. (219.11,173.11) .. controls (221.39,172.5) and (222.83,173.33) .. (223.44,175.61) .. controls (224.05,177.89) and (225.49,178.72) .. (227.77,178.11) .. controls (230.05,177.5) and (231.49,178.33) .. (232.1,180.61) .. controls (232.71,182.89) and (234.15,183.72) .. (236.43,183.11) .. controls (238.71,182.5) and (240.15,183.33) .. (240.76,185.61) .. controls (241.37,187.89) and (242.81,188.72) .. (245.09,188.11) .. controls (247.37,187.5) and (248.81,188.33) .. (249.42,190.61) .. controls (250.03,192.88) and (251.48,193.71) .. (253.75,193.1) .. controls (256.03,192.49) and (257.47,193.32) .. (258.08,195.6) .. controls (258.69,197.88) and (260.13,198.71) .. (262.41,198.1) .. controls (264.69,197.49) and (266.13,198.32) .. (266.74,200.6) .. controls (267.35,202.88) and (268.8,203.71) .. (271.08,203.1) .. controls (273.36,202.49) and (274.8,203.32) .. (275.41,205.6) .. controls (276.02,207.88) and (277.46,208.71) .. (279.74,208.1) .. controls (282.02,207.49) and (283.46,208.32) .. (284.07,210.6) .. controls (284.68,212.87) and (286.13,213.7) .. (288.4,213.09) .. controls (290.68,212.48) and (292.12,213.31) .. (292.73,215.59) .. controls (293.34,217.87) and (294.78,218.7) .. (297.06,218.09) .. controls (299.34,217.48) and (300.78,218.31) .. (301.39,220.59) .. controls (302,222.87) and (303.44,223.7) .. (305.72,223.09) .. controls (308,222.48) and (309.44,223.31) .. (310.05,225.59) .. controls (310.66,227.87) and (312.1,228.7) .. (314.38,228.09) .. controls (316.66,227.48) and (318.11,228.31) .. (318.72,230.59) .. controls (319.33,232.86) and (320.78,233.69) .. (323.05,233.08) .. controls (325.33,232.47) and (326.77,233.3) .. (327.38,235.58) .. controls (327.99,237.86) and (329.43,238.69) .. (331.71,238.08) .. controls (333.99,237.47) and (335.43,238.3) .. (336.04,240.58) .. controls (336.65,242.86) and (338.09,243.69) .. (340.37,243.08) .. controls (342.65,242.47) and (344.09,243.3) .. (344.7,245.58) .. controls (345.31,247.86) and (346.75,248.69) .. (349.03,248.08) .. controls (351.31,247.47) and (352.75,248.3) .. (353.36,250.58) .. controls (353.97,252.85) and (355.42,253.68) .. (357.69,253.07) .. controls (359.97,252.46) and (361.41,253.29) .. (362.02,255.57) .. controls (362.63,257.85) and (364.07,258.68) .. (366.35,258.07) -- (368.03,259.04) -- (368.03,259.04) ;
\path  [shading=_lwdhe5o8o,_bgoua249b,opacity=0.6] (149.81,133.13) -- (261.76,311.97) -- (152.08,259.04) -- cycle ; 
 \draw   (149.81,133.13) -- (261.76,311.97) -- (152.08,259.04) -- cycle ; 

\path  [shading=_1gw1wh7zt,_flungya80,opacity=0.6] (261.76,311.67) -- (369,132.67) -- (368.03,259.04) -- cycle ; 
 \draw   (261.76,311.67) -- (369,132.67) -- (368.03,259.04) -- cycle ; 

\path  [shading=_otc30b6nt,_ff01gojv6,opacity=0.6] (259.41,78) -- (369,133.13) -- (149.81,133.13) -- cycle ; 
 \draw   (259.41,78) -- (369,133.13) -- (149.81,133.13) -- cycle ; 

\draw (130,125) node [anchor=north west][inner sep=0.75pt]   [align=left] {$12$};
\draw (252,58) node [anchor=north west][inner sep=0.75pt]   [align=left] {$23$};
\draw (190,114) node [anchor=north west][inner sep=0.75pt]   [align=left] {\textcolor{magenta}{$F_{123|4}$}};
\draw (176,138) node [anchor=north west][inner sep=0.75pt]   [align=left] {\textcolor{gray}{$F_{12|34}$}};
\draw (200,196) node [anchor=north west][inner sep=0.75pt]   [align=left] {\textcolor{gray}{$F_{13|24}$}};
\draw (375,125) node [anchor=north west][inner sep=0.75pt]   [align=left] {$13$};
\draw (375,250) node [anchor=north west][inner sep=0.75pt]   [align=left] {$34$};
\draw (268,144) node [anchor=north west][inner sep=0.75pt]   [align=left] {\textcolor{magenta}{$F_{234|1}$}};
\draw (252,318) node [anchor=north west][inner sep=0.75pt]   [align=left] {$14$};
\draw (130,250) node [anchor=north west][inner sep=0.75pt]   [align=left] {$24$};
\draw (206,272) node [anchor=north west][inner sep=0.75pt]   [align=left] {\textcolor{magenta}{$F_{124|3}$}};
\draw (262,228) node [anchor=north west][inner sep=0.75pt]   [align=left] {\textcolor{gray}{$F_{14|23}$}};
\draw (283,271) node [anchor=north west][inner sep=0.75pt]   [align=left] {\textcolor{magenta}{$F_{134|2}$}};
\end{tikzpicture}
\caption{The complex  \(\Delta(\scrA_3)= \Delta(\frs\frl_4)\).}
\end{figure}
\end{center}

It follows from the lattice structure and  \eqref{eq:dimSRR} that 
\[
\dim \Bbbk[\Delta(\frs\frl_4)] \ = \  \max\{ \# F_\pi \mid \pi \in L(\scrA)^\circ \} = 3 = \operatorname{rank} \frs\frl_4.
\] 
\end{example}

The main result in this section is the following.

\setcounter{table}{0}
\begin{theorem}
\label{thm:SRdim}
If \( \frg \) is a classical simple Lie algebra of rank \( n\geq 2\), and \( \Delta(\frg) \) is the nerve complex of its Coxeter arrangement, then the dimension of the corresponding Stanley-Reisner ring \( \Bbbk[\Delta(\frg)] \) is displayed in the following table.
\begin{table}[hpt]
\begin{center}
\begin{tabular}{|c|c|c|}\hline
Type & \( \dim \Bbbk[\Delta(\frg)]\) & \( \operatorname{codim} I_{\Delta(\frg)} \) \\ \hline \hline
\( A_n\) & \( \binom{n}{2}\) & \( n \) \\
\( B_n\) & \( (n-1)^2 \)& \( 2n -1\) \\
\( C_n\) & \( (n-1)^2\) & \( 2n -1\) \\
\( D_n\) & \( 2\binom{n-1}{2} \)   & \( 2(n -1)\)\\ \hline
\end{tabular}
\end{center}
\caption{Dimension and codimension of Stanley-Reisner ideals}
\label{dim-SRn}
\end{table}%
\end{theorem}

\subsection{Signed graphs, partition complexes and the proof of Theorem \ref{thm:SRdim}}

Using   Remark \ref{rem:SR}, it suffices to determine 
the maximal number of vertices in a facet of \( \Delta( \frg)\), i.e.  \( \max \{ \# F_\pi \mid \pi \in L(\scrA(\frg))^\circ\}.
\) As a result, the proof follows from  the well-known structure of the intersection lattice \( L(\scrA(\frg))\) of the Coxeter arrangement, which has various description in the literature. Here we follow \cite{MR1406454}, except for a minor difference in notation, explained below.

\begin{notation}
\label{not:ABCD}
Given \( n\geq 2 \)  consider the following arrangements of hyperplanes in \( \bbr^n\), and note that \( \scrA_{n-1} \subset \scrD_n \subset \scrB_n\).
\begin{align}
 \scrA_{n-1} & := \{ x_i = x_j \mid 1 \leq i< j \leq n \} \notag \\
 \scrB_n  & := \{ x_i = \pm x_j \mid 1 \leq i< j \leq n \} \cup \{ x_i = 0 \mid 1 \leq i<\leq n \}\\
 \scrD_n & := \{ x_i = \pm x_j \mid 1 \leq i< j \leq n \}. \notag
\end{align}
In \cite{MR1406454} the authors use \( \cala_n \) to denote what we call \( \scrA_{n-1} \) above.
\end{notation}

Now, let us recall that a signed graph \( \Gamma \) on a (finite) vertex set \( V(\Gamma) \) can have the following types of edges:
\begin{enumerate}[1.]
\item A \emph{positive edge} \( ij^+\) between two distinct vertices \( i, j \in V(\Gamma) \). 
\item A \emph{negative edge} \( ij^-\) between two distinct vertices \( i, j \in V(\Gamma) \). 
\item A \emph{loop} (or a half-edge) \( i^h\)  based on the vertex \( i  \in V(\Gamma) \). Such a loop is seen as a having a negative sign.
\end{enumerate}

\begin{example}
\label{exmp:intuition}
The intuition  comes from the complete signed graph  \( G_{\scrB_n} \) on \([n]\),  associated to the arrangement \( \scrB_n\). Here, the edge \( ij^+ \) corresponds to the hyperplane \( x_i = x_j\) and \( ij^- \) corresponds to the hyperplane \( x_i = -x_j\), when \( i\neq j \), whereas the loop \( i^h\) corresponds to the hyperplane \( x_i = 0 \).
\end{example}

\begin{definition}\cite{MR1406454}
\label{def:complete}
Let  \( \scrA \subseteq \scrB_n\) be a subarrangement on   the vertex set \( [n] \).

\begin{enumerate}[i.]
\item  Define the \emph{associated signed graph} \( G_\scrA\) by saying that an edge of \( G_{\scrB_n} \) is in  \( G_\scrA\)   if and only if the corresponding hyperplane is in \( \scrA\).
\item Now, restrict \( \scrA \) to be one of the three arrangements \( \scrA_n \subset \scrD_n\subset \scrB_n\),
and consider subsets \( V, W \subset [n] \), with \( V \cap W = \emptyset \) and \( V \cup W \neq \emptyset \).
\begin{enumerate}
\item The \emph{complete balanced graph} \( K^b_{V,W} \)  consists of all positive edges between vertices of \( V \), all positive edges between vertices of \( W\), and all negative edges from vertices of  \( V\) to vertices of \( W\).
\item The \emph{complete unbalanced graph} \( K^u_{V} \) consists of all positive and negative edges between vertices of \( V \). When \( \scrA = \scrB_n\), the edges of \(K^u_V\) also include the loops based on each vertex of \( V \). 
\end{enumerate}
\end{enumerate}
\end{definition}

\begin{figure}[ht]
\centering
\begin{subfigure}{.44\textwidth}
  \centering


\tikzset{every picture/.style={line width=0.75pt}} 

\begin{tikzpicture}[x=0.44pt,y=0.44pt,yscale=-1,xscale=1]

\draw   (226.5,275) .. controls (226.5,272.24) and (224.26,270) .. (221.5,270) .. controls (218.74,270) and (216.5,272.24) .. (216.5,275) .. controls (216.5,277.76) and (218.74,280) .. (221.5,280) .. controls (224.26,280) and (226.5,277.76) .. (226.5,275) -- cycle ;
\draw   (512.5,272) .. controls (512.5,269.24) and (510.26,267) .. (507.5,267) .. controls (504.74,267) and (502.5,269.24) .. (502.5,272) .. controls (502.5,274.76) and (504.74,277) .. (507.5,277) .. controls (510.26,277) and (512.5,274.76) .. (512.5,272) -- cycle ;
\draw   (381.5,197) .. controls (381.5,194.24) and (379.26,192) .. (376.5,192) .. controls (373.74,192) and (371.5,194.24) .. (371.5,197) .. controls (371.5,199.76) and (373.74,202) .. (376.5,202) .. controls (379.26,202) and (381.5,199.76) .. (381.5,197) -- cycle ;
\draw   (281.55,405) .. controls (281.55,402.24) and (280.19,400) .. (278.52,400) .. controls (276.85,400) and (275.5,402.24) .. (275.5,405) .. controls (275.5,407.76) and (276.85,410) .. (278.52,410) .. controls (280.19,410) and (281.55,407.76) .. (281.55,405) -- cycle ;
\draw   (454.5,402) .. controls (454.5,399.24) and (453.15,397) .. (451.48,397) .. controls (449.81,397) and (448.45,399.24) .. (448.45,402) .. controls (448.45,404.76) and (449.81,407) .. (451.48,407) .. controls (453.15,407) and (454.5,404.76) .. (454.5,402) -- cycle ;
\draw [line width=1.5]    (221.5,275) -- (376.5,197) ;
\draw [line width=1.5]    (376.5,197) -- (507.5,272) ;
\draw [line width=1.5]    (221.5,275) -- (507.5,272) ;
\draw [line width=1.5]    (278.52,405) -- (451.48,402) ;
\draw [color={rgb, 255:red, 208; green, 2; blue, 27 }  ,draw opacity=1 ]   (278.52,405) -- (221.5,275) ;
\draw [color={rgb, 255:red, 208; green, 2; blue, 27 }  ,draw opacity=0.6 ] [dash pattern={on 3pt off 2pt}]  (278.52,405) -- (376.5,197) ;
\draw [color={rgb, 255:red, 208; green, 2; blue, 27 }  ,draw opacity=1 ]   (278.52,405) -- (507.5,272) ;
\draw [color={rgb, 255:red, 208; green, 2; blue, 27 }  ,draw opacity=0.6 ] [dash pattern={on 3pt off 2pt}]  (451.48,402) -- (376.5,197) ;
\draw [color={rgb, 255:red, 208; green, 2; blue, 27 }  ,draw opacity=1 ]   (451.48,402) -- (507.5,272) ;
\draw [color={rgb, 255:red, 208; green, 2; blue, 27 }  ,draw opacity=1 ]   (362.5,354) -- (221.5,275) ;
\draw [color={rgb, 255:red, 208; green, 2; blue, 27 }  ,draw opacity=1 ]   (448.45,402) -- (369.5,357) ;

\draw (190,268.4) node [anchor=north west][inner sep=0.75pt]    {\footnotesize $1$};
\draw (371,168.4) node [anchor=north west][inner sep=0.75pt]    {\footnotesize$3$};
\draw (528,265.4) node [anchor=north west][inner sep=0.75pt]    {\footnotesize$2$};
\draw (251,396.4) node [anchor=north west][inner sep=0.75pt]    {\footnotesize$4$};
\draw (464,392.4) node [anchor=north west][inner sep=0.75pt]    {\footnotesize$5$};
\draw (265,215) node [anchor=north west][inner sep=0.75pt]    {\footnotesize$13^+$};
\draw (445,215) node [anchor=north west][inner sep=0.75pt]    {\footnotesize$23^+$};
\draw (354,277.4) node [anchor=north west][inner sep=0.75pt]    {\footnotesize $12^+$};
\draw (356,408.4) node [anchor=north west][inner sep=0.75pt]    {\footnotesize $45^+$};
\draw (206,320) node [anchor=north west][inner sep=0.75pt]    {\footnotesize $\textcolor[rgb]{0.82,0.01,0.11}{14^-}$};
\draw (275,338) node [anchor=north west][inner sep=0.75pt]    {\footnotesize $\textcolor[rgb]{0.82,0.01,0.11}{34^-}$};
\draw (485,320) node [anchor=north west][inner sep=0.75pt]    {\footnotesize $\textcolor[rgb]{0.82,0.01,0.11}{25^-}$};
\draw (436,338) node [anchor=north west][inner sep=0.75pt]    {\footnotesize $\textcolor[rgb]{0.82,0.01,0.11}{35^-}$};
\draw (318,379.4) node [anchor=north west][inner sep=0.75pt]    {\footnotesize $\textcolor[rgb]{0.82,0.01,0.11}{24^-}$};
\draw (380,379.4) node [anchor=north west][inner sep=0.75pt]    {\footnotesize $\textcolor[rgb]{0.82,0.01,0.11}{15^-}$};
\end{tikzpicture}
 \caption*{\( K^b_{\{ 1,2,3\}, \,  \{ 4,5\} }\)}
  \label{fig:sub1}
\end{subfigure}%
\begin{subfigure}{.56\textwidth}
  \centering
%


\tikzset{every picture/.style={line width=0.75pt}} 

\begin{tikzpicture}[x=0.28pt,y=0.29pt,yscale=-1,xscale=1]

\draw   (458.13,247.12) .. controls (458.13,242.76) and (461.2,239.23) .. (464.99,239.23) .. controls (468.78,239.23) and (471.85,242.76) .. (471.85,247.12) .. controls (471.85,251.47) and (468.78,255) .. (464.99,255) .. controls (461.2,255) and (458.13,251.47) .. (458.13,247.12) -- cycle ;
\draw  [dash pattern={on 3pt off 2pt}, draw opacity=0.4]  (242.35,488.38) .. controls (333.51,449.78) and (577.11,444.95) .. (684.4,481.95) ;
\draw [color={rgb, 255:red, 208; green, 2; blue, 27 }  ,draw opacity=0.4 ] [dash pattern={on 3pt off 3pt}]  (242.35,488.38) .. controls (331.89,514.12) and (583.57,512.51) .. (684.4,481.95) ;
\draw    (242.35,488.38) .. controls (270.59,407.96) and (349.64,293.76) .. (464.99,247.12) ;
\draw [color={rgb, 255:red, 208; green, 2; blue, 27 }  ,draw opacity=1 ]   (242.35,488.38) .. controls (333.51,424.05) and (433.53,311.46) .. (464.99,247.12) ;
\draw [color={rgb, 255:red, 0; green, 0; blue, 0 }  ,draw opacity=1 ]   (684.4,481.95) .. controls (661,404.74) and (551.3,274.46) .. (464.99,247.12) ;
\draw [color={rgb, 255:red, 208; green, 2; blue, 27 }  ,draw opacity=1 ]   (684.4,481.95) .. controls (620.67,441.74) and (498.06,325.93) .. (464.99,247.12) ;
\draw    (516.61,629.92) .. controls (419.01,617.06) and (296.4,562.37) .. (242.35,488.38) ;
\draw [color={rgb, 255:red, 208; green, 2; blue, 27 }  ,draw opacity=1 ]   (516.61,629.92) .. controls (441.6,547.89) and (362.54,523.77) .. (242.35,488.38) ;
\draw [color={rgb, 255:red, 208; green, 2; blue, 27 }  ,draw opacity=1 ]   (684.4,481.95) .. controls (585.18,530.2) and (559.37,555.94) .. (515,628.31) ;
\draw    (684.4,481.95) .. controls (664.23,539.85) and (573.89,612.23) .. (516.61,629.92) ;
\draw [color={rgb, 255:red, 208; green, 2; blue, 27 }  ,draw opacity=1 ]   (515,628.31) .. controls (473.33,534) and (436.72,397) .. (464.99,247.12) ;
\draw [color={rgb, 255:red, 0; green, 0; blue, 0 }  ,draw opacity=1 ]   (515,628.31) .. controls (538.39,475.51) and (481.93,308.24) .. (464.99,247.12) ;
\draw   (677.54,481.95) .. controls (677.54,477.59) and (680.61,474.06) .. (684.4,474.06) .. controls (688.18,474.06) and (691.25,477.59) .. (691.25,481.95) .. controls (691.25,486.3) and (688.18,489.83) .. (684.4,489.83) .. controls (680.61,489.83) and (677.54,486.3) .. (677.54,481.95) -- cycle ;
\draw   (508.14,628.31) .. controls (508.14,623.96) and (511.21,620.43) .. (515,620.43) .. controls (518.79,620.43) and (521.86,623.96) .. (521.86,628.31) .. controls (521.86,632.67) and (518.79,636.2) .. (515,636.2) .. controls (511.21,636.2) and (508.14,632.67) .. (508.14,628.31) -- cycle ;
\draw   (235.5,488.38) .. controls (235.5,484.03) and (238.57,480.5) .. (242.35,480.5) .. controls (246.14,480.5) and (249.21,484.03) .. (249.21,488.38) .. controls (249.21,492.74) and (246.14,496.27) .. (242.35,496.27) .. controls (238.57,496.27) and (235.5,492.74) .. (235.5,488.38) -- cycle ;
\draw  [draw opacity=0] (454.68,246.14) .. controls (442.66,241.14) and (433.94,226.16) .. (433.94,208.46) .. controls (433.94,186.67) and (447.17,169) .. (463.48,169) .. controls (479.79,169) and (493.02,186.67) .. (493.02,208.46) .. controls (493.02,224.59) and (485.78,238.46) .. (475.4,244.58) -- (463.48,208.46) -- cycle ; \draw  [color={rgb, 255:red, 208; green, 2; blue, 27 }  ,draw opacity=1 ] (454.68,246.14) .. controls (442.66,241.14) and (433.94,226.16) .. (433.94,208.46) .. controls (433.94,186.67) and (447.17,169) .. (463.48,169) .. controls (479.79,169) and (493.02,186.67) .. (493.02,208.46) .. controls (493.02,224.59) and (485.78,238.46) .. (475.4,244.58) ;  
\draw  [draw opacity=0] (526.42,633.55) .. controls (538.36,638.69) and (546.88,653.77) .. (546.64,671.47) .. controls (546.34,693.26) and (532.87,710.77) .. (516.55,710.57) .. controls (500.24,710.38) and (487.26,692.56) .. (487.56,670.77) .. controls (487.78,654.64) and (495.21,640.86) .. (505.67,634.86) -- (517.1,671.12) -- cycle ; \draw   [color={rgb, 255:red, 208; green, 2; blue, 27 }  ,draw opacity=1 ]  (526.42,633.55) .. controls (538.36,638.69) and (546.88,653.77) .. (546.64,671.47) .. controls (546.34,693.26) and (532.87,710.77) .. (516.55,710.57) .. controls (500.24,710.38) and (487.26,692.56) .. (487.56,670.77) .. controls (487.78,654.64) and (495.21,640.86) .. (505.67,634.86) ;  
\draw  [draw opacity=0] (688.39,469.02) .. controls (694.33,458.48) and (710.02,450.98) .. (728.4,451.03) .. controls (751.87,451.09) and (770.85,463.42) .. (770.81,478.58) .. controls (770.76,493.73) and (751.7,505.96) .. (728.24,505.9) .. controls (711.74,505.86) and (697.46,499.75) .. (690.45,490.86) -- (728.32,478.46) -- cycle ; \draw  [color={rgb, 255:red, 208; green, 2; blue, 27 }  ,draw opacity=1 ] (688.39,469.02) .. controls (694.33,458.48) and (710.02,450.98) .. (728.4,451.03) .. controls (751.87,451.09) and (770.85,463.42) .. (770.81,478.58) .. controls (770.76,493.73) and (751.7,505.96) .. (728.24,505.9) .. controls (711.74,505.86) and (697.46,499.75) .. (690.45,490.86) ;  
\draw  [draw opacity=0] (235.51,501.26) .. controls (229.76,511.9) and (214.22,519.65) .. (195.83,519.9) .. controls (172.37,520.22) and (153.16,508.19) .. (152.93,493.04) .. controls (152.69,477.88) and (171.52,465.34) .. (194.98,465.03) .. controls (211.48,464.81) and (225.87,470.69) .. (233.04,479.47) -- (195.41,492.46) -- cycle ; \draw  [color={rgb, 255:red, 208; green, 2; blue, 27 }  ,draw opacity=1 ] (235.51,501.26) .. controls (229.76,511.9) and (214.22,519.65) .. (195.83,519.9) .. controls (172.37,520.22) and (153.16,508.19) .. (152.93,493.04) .. controls (152.69,477.88) and (171.52,465.34) .. (194.98,465.03) .. controls (211.48,464.81) and (225.87,470.69) .. (233.04,479.47) ;  

\draw (200,475) node [anchor=north west][inner sep=0.75pt]    {\footnotesize $1$};
\draw (504,650) node [anchor=north west][inner sep=0.75pt]    {\footnotesize $2$};
\draw (700,466) node [anchor=north west][inner sep=0.75pt]    {\footnotesize $3$};
\draw (450,205) node [anchor=north west][inner sep=0.75pt]    {\footnotesize $4$};
\draw (280,315) node [anchor=north west][inner sep=0.75pt]    {\footnotesize $14^+$};
\draw (305,370) node [anchor=north west][inner sep=0.75pt]  [color={rgb, 255:red, 208; green, 2; blue, 27 }  ,opacity=1 ]  {\footnotesize $14^-$};
\draw (616,344.37) node [anchor=north west][inner sep=0.75pt]    {\footnotesize $34+$};
\draw (536,330) node [anchor=north west][inner sep=0.75pt]  [color={rgb, 255:red, 208; green, 2; blue, 27 }  ,opacity=1 ]  {\footnotesize $34-$};
\draw (510,407.64) node [anchor=north west][inner sep=0.75pt]    {\footnotesize $24+$};
\draw (400,380) node [anchor=north west][inner sep=0.75pt]  [color={rgb, 255:red, 208; green, 2; blue, 27 }  ,opacity=1 ]  {\footnotesize $24^-$};
\draw (570,546.64) node [anchor=north west][inner sep=0.75pt]  [color={rgb, 255:red, 208; green, 2; blue, 27 }  ,opacity=1 ]  {\footnotesize $23^-$};
\draw (635,546.64) node [anchor=north west][inner sep=0.75pt]  [color={rgb, 255:red, 0; green, 0; blue, 0 }  ,opacity=1 ]  {\footnotesize $23^+$};
\draw (371.81,550) node [anchor=north west][inner sep=0.75pt]  [color={rgb, 255:red, 208; green, 2; blue, 27 }  ,opacity=1 ]  {\footnotesize $12^-$};
\draw (344.93,600) node [anchor=north west][inner sep=0.75pt]  [color={rgb, 255:red, 0; green, 0; blue, 0 }  ,opacity=1 ]  {\footnotesize $12^+$};
\draw (383.65,476) node [anchor=north west][inner sep=0.75pt]  [color={rgb, 255:red, 208; green, 2; blue, 27 }  ,opacity=1 ]  {\footnotesize $13^-$};
\draw (350,430) node [anchor=north west][inner sep=0.75pt]  [color={rgb, 255:red, 0; green, 0; blue, 0 }  ,opacity=1 ]  {\footnotesize $13^+$};
\draw (500,170) node [anchor=north west][inner sep=0.75pt]    {\footnotesize $\textcolor[rgb]{0.82,0.01,0.11}{ 4^h}$};
\draw (730,510) node [anchor=north west][inner sep=0.75pt]    {\footnotesize $\textcolor[rgb]{0.82,0.01,0.11}{ 3^h}$};
\draw (444,670) node [anchor=north west][inner sep=0.75pt]    {\footnotesize $\textcolor[rgb]{0.82,0.01,0.11}{  2^h}$};
\draw (180,425) node [anchor=north west][inner sep=0.75pt]    {\footnotesize $\textcolor[rgb]{0.82,0.01,0.11}{  1^h}$};
\end{tikzpicture}
\caption*{ \( K^u_{\{1, 2, 3, 4\}}\)}
  \label{fig:sub2}
\end{subfigure}
\caption{Complete balanced and unbalanced signed graphs}
\label{fig:test}
\end{figure}
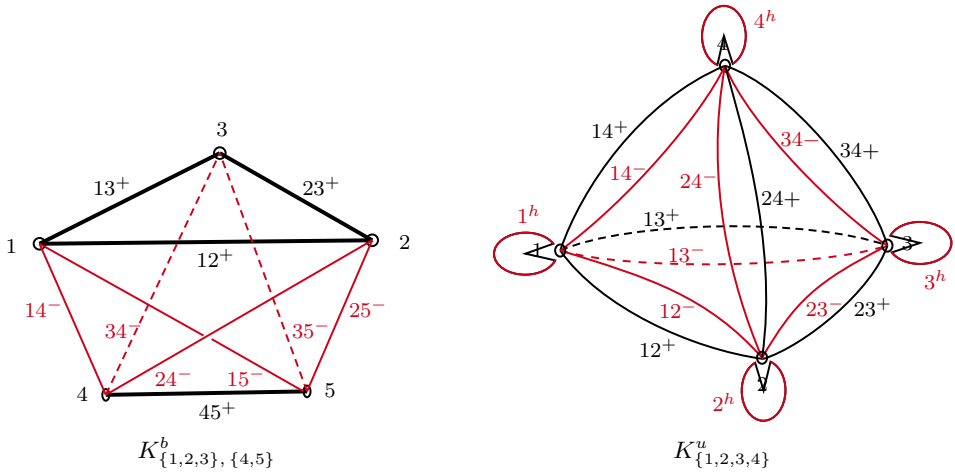

%
%

\begin{theorem}[\cite{MR676405}, \cite{MR1406454}]
\label{thm:class}
Let \( \scrA = \scrA_{n-1}, \scrB_n\) or \(\scrD_n \).  The lattice \( L(\scrA)\) is isomorphic to the lattice of subgraphs \( G \subset G_\scrA\) such that
\begin{enumerate}[1.]
\item every component of \( G \) is complete balanced or complete unbalanced, and
\item there is at most one unbalanced component.
\end{enumerate}
The graphs are partially ordered by inclusion of their edge sets. 
\end{theorem}

\begin{proof}[Proof of Theorem \ref{thm:SRdim}]
We need to determine the dimension of \( \Delta(\scrA) \), with \( \scrA\) as in Theorem \ref{thm:class}. By definition,   facets of \( \Delta(\scrA)\) come from  maximal elements in \( L(\scrA)^\circ \) which, by Theorem \ref{thm:class},  correspond to proper subgraphs \( G \subsetneq G_\scrA\)   satisfying the two conditions in the theorem, and such that 
\( V(G) = V(G_\scrA) = [n] \).
\bigskip

\noindent{\bf The case \( \scrA_n\).}\ \ The facets in this case correspond to non-trivial partitions \( \pi :  P_1 \amalg \cdots \amalg P_r = [n+1] \), with 
\( r \geq 2 \) and \( P_j \neq \emptyset\) for all \( j=1, \ldots, r.\)  Since there are only positive edges in this case, one has
\( G_\pi = K^b_{P_1, \emptyset}\amalg \cdots \amalg K^b_{P_r, \emptyset}\), and the edge set is given by 
\(
E(G_\pi) = \Lwedge^2 P_1 \amalg \cdots \amalg \Lwedge^2 P_r.
\)
Note that when \( \# P_i = 1 \) then \( \Lwedge^2P_i = \emptyset \). Since \( r \geq 2 \), it follows that
\[
\# F_\pi = \sum_{\ell=1}^r \binom{\# P_\ell}{2} \leq \binom{\left\{ \sum_{\ell=1}^r \#P_\ell \right\} + 1 -r}{2}\ = \ \binom{n+2-r}{2} \leq \binom{n}{2}.
\]
Note that equality is attained when \( \pi \) is a partition of the form \( \{ i \} \amalg ( [n+1]-\{i\} )\). Therefore, 
\( \dim(\Bbbk[\Delta(\scrA_n)] = \binom{n}{2} \) and \( \operatorname{codim}(I_{\Delta(\scrA_n)}) = \# R^+ - \binom{n}{2} = 
\binom{n+1}{2}-\binom{n}{2} = n \), according to Remark \ref{rem:SR}. This proves the  case \( \scrA_n\) of the theorem.
\medskip

\noindent{\bf The case \( \scrD_n\).}\ \ In this case, the subgraph \( G_\pi\)    corresponding to a facet \( \pi\) has the form
\( G_\pi = \left( K^b_{V_1,W_1} \amalg \cdots \amalg K^b_{V_r,W_r} \right) \amalg K^u_{S_{r+1}} \)
where
\( \pi :  \left( V_1 \amalg W_1 \right) \amalg \cdots \amalg \left( V_r \amalg W_r \right) \amalg S_{r+1}
\)
is a partition of \( [n] \), satisfying \( V_i \cup W_i \neq \emptyset,\ i = 1, \ldots, r \), whereas \( S_{r+1} \) may be empty. 
Denote \( \# V_i = v_i, \ \# W_i = w_i \) and \( \# S_{r+1} = s \), so that \(\sum_{i=1}^r (v_i+w_i) \ + \ s \ = \ n.\)

It is easy to see that \( \# E(K^b_{V,W} ) = \binom{\#V\,  +\, \#W}{2} \) and \( \# E(K^u_S) = 2 \binom{\#S}{2} \), whenever  \(\scrA =  \scrD_n\). 
Therefore, if we denote \( m_\ell = v_\ell + w_\ell\), then
\[
\# E(G_\pi) \ = \ \sum_{\ell=1}^r \binom{m_\ell}{2} \ + \ 2 \binom{s}{2}.
\]
Since \( \sum_{\ell=1}^r m_\ell = n-s \) with \( 1 \leq m_\ell \leq n\),  \( 0 \leq s \leq n-1\) and \( r \geq 1 \), one has
\[
\# E(G_\pi) \leq \binom{n-s+1-r}{2} + 2 \binom{s}{2} \leq 2 \binom{n-1}{2},
\]
and equality is attained with  subgraphs of the form \( G_\pi =  K^b_{\{ i \}, \emptyset} \amalg K^u_{[n] - \{ i\}} \), for any \( i=1, \ldots, n\). This proves the case \( \scrD_n \) of the theorem.
\medskip

\noindent{\bf The cases \( \scrB_n\) and \( \scrC_n\).}\ \  Just like in the previous case, facets of \( \Delta(\scrB_n) \) correspond to subgraphs of the type
\(
 G = \left( K^b_{V_1,W_1} \amalg \cdots \amalg K^b_{V_r,W_r} \right) \amalg K^u_{S_{r+1}}.
\)
However, in this case \( \# E(K^u_S )= 2\binom{s}{2} + s  = s^2,\) and hence
\begin{align*}
\# E(G_\pi) & = \# E(G_\pi) \ = \ \sum_{\ell=1}^r \binom{m_\ell}{2} \ + s^2
\\
& \leq \binom{n-s+1-r}{2} + s^2 \ \leq \ \binom{n-s}{2} + s^2 \leq (n-1)^2 .
\end{align*}
Note that equality is attained with subgraphs
of the form \( G_\pi =  K^b_{\{ i \}, \emptyset} \amalg K^u_{[n] - \{ i\}} \), for any \( i=1, \ldots, n\). This proves the cases \( \scrB_n \) and \( \scrC_n\) of the theorem, since the lattices \( L(\scrB_n)\) and \(L(\scrC_n)\) are isomorphic.
 \end{proof}



\section{LULU is equidimensional}
\label{sec:equidim}

We have defined   ideals  \( \cali_\calb \) and  \( I_{\Delta(\frg)}\)  in the polynomial ring \( \Bbbk[x_1, \ldots, x_D]\), where \( x_i \) is a variable associated to a positive root \( \alpha_i \) in \( R^+ = \{ \alpha_1, \ldots, \alpha_D\} \); see \eqref{eq:iso}. The   ideal  \( \cali_\calb \) was introduced in \S \ref{subsec:pb} and defines the pull-back of a determinantal scheme under the morphism \( \varphi\),  while  \( I_{\Delta(\frg)}\)  is the Stanley-Reisner ideal  associated to the Coxeter arrangement \( \scrA(\frg)\); see Definition \ref{def:nerve}.
The first objective of this section is to establish the following  relationship between these two ideals. 

\begin{theorem}
\label{thm:inclusion}
Let \( \frg \) be a simple Lie algebra of rank \ \( n \) over \( \Bbbk \). Using the notation above, we have an inclusion of ideals \ 
\( I_{\Delta(\frg)} \ \subseteq \  \cali_\calb \subset \ \Bbbk[\mbx] . \)
\end{theorem}
Once this result is shown, we proceed to prove Theorem \ref{thm:main}.

\subsection{The relationship between \( \cali_\calb \) and \( I_{\Delta(\frg)}\)}

Consider the composition
\begin{equation}
\label{eq:F}
\begin{tikzcd}
\bba^D \ar[rrrr, rounded corners=8pt, to path={-- ([yshift=-3.5ex]\tikztostart.south)-| (\tikztotarget)}] 
\ar[r,"\Psi", "\cong"'] & \bbu \arrow[r, hook] 
& 
G \arrow[r, "\mathsf{Ad}"] \ar[d, phantom, "\scriptstyle F", pos=1] 
 & 
GL(\frg) \arrow[r, "\mathsf{res}_{\frl}"]    &
\HomSch{\frl}{\frg} \\
& & \phantom{.} & &  
\end{tikzcd}
\end{equation}
so that \( \varphi\circ \Psi^{-1} = \mathsf{proj}_{\frh} \circ F \colon \bbu \to \HomSch{\frl}{\frh} \); see \eqref{eq:varphi}.  As explained in Remark \ref{rem:Hom}, the morphism \( F \) corresponds to an element \( \check{F} \in  \frl^*\otimes \frg \otimes \Bbbk[\mbx]\) that can be written as a  sum of homogeneous components \( \check{F} = \check{F}_0 + \cdots + \check{F}_m  \), with \( \check{F}_d \in  \frl^*\otimes \frg \otimes \Bbbk[\mbx]_{(d)}\), and where \( \Bbbk[\mbx]_{(d)} \) denotes the  space of homogeneous polynomials of degree \( d\).

If \( I \in \frl^*\otimes \frg \) corresponds to the inclusion \( \frl \hookrightarrow \frg\), then \( \check{F}_0 = I \). To determine 
\( \check{F}_1 \), first observe that the differential 
\( dF_{|0} \colon T_0\bba^D \to T_0 \HomSch{\frl}{\frg} = \Hom{\Bbbk}{\frl}{\frg}\)  sends \( \frac{\partial}{\partial x_\beta} \)
to \( \ad{X_\beta}_{|\frl}\) for all \( \beta \in R^+\); see \eqref{eq:diffPsi}. Therefore, 
\begin{equation}
\label{eq:F1}
\check{F}_1 = \sum_{\beta\in R^+}\, \ad{X_\beta}_{|\frl} \otimes x_\beta.
\end{equation}

In particular, for any \( \alpha \in R^+\) one has 
\begin{align*}
\check{F}_1(Y_\alpha) & =
\sum_{\beta \in R^+}\, [X_\beta, Y_\alpha] \otimes x_\beta :=
\sum_{\beta \in R^+}\, [X_\beta, X_{-\alpha}] \otimes x_\beta \\
& =\calb
[X_\alpha, X_{-\alpha}] \otimes x_\alpha \ + \ 
\sum_{\substack{\beta - \alpha \in R \\ \beta \neq \alpha}} \ N_{\beta, - \alpha}\ X_{\beta -\alpha}\otimes x_\beta\\
& =
H_\alpha \otimes x_\alpha \ + \ 
\sum_{\substack{\beta - \alpha \in R \\ \beta \neq \alpha}} \ N_{\beta, - \alpha}\ X_{\beta -\alpha}\otimes x_\beta\\
\end{align*}
where \( N_{\beta, -\alpha} \) are the structural constants of \( \frg\) and we denote \( Y_{ \alpha} = X_{-\alpha} \) interchangeably; see \cite{MR3616493}.
Hence, the homogeneous component of degree \( 1 \) 
of \( \varphi_{\Psi(\mbx)} \),   denoted by \( \varphi^{(1)}_\mbx\), is completely determined by the fact that for all \( \alpha \in R^+ \) one has
\begin{equation}
\label{eq:varphi1}
\varphi_{\mbx}^{(1)}(Y_\alpha) =   \proj{\frh}{\check{F}_1(Y_\alpha)}   = H_\alpha \otimes 
x_\alpha .
\end{equation}

Define 
\begin{equation*}
\cals_n(\frg) = \left\{ J = (j_1< \cdots < j_n) \in \Lwedge^n[D] \mid \{ \alpha_{j_1}, \ldots, \alpha_{j_n} \} \text{ is a linear basis for } \frh^* \right\} . 
\end{equation*}
Given an element \( J \in \cals_n(\frg) \), fix the corresponding basis \( \{ H_{\alpha_{j_1}}, \ldots, H_{\alpha_{j_n}} \} \) for \( \frh\).

Now, recall the set of monomial generators of \(\cali_{\calb}\) explained in \eqref{eq:caliB}, whose description started by fixing a basis for \( \frh\). If the fixed basis is the one above, consider
\( q_J(\mbx) = \det \left( p_{j_r, s}(\mbx)  \right) \) as in \eqref{eq:qI}, where \(  p_{j_r, s}(\mbx) \) comes from writing
\begin{equation}
\label{eq:varphi2}
\varphi_{\Psi(\mbx)} = \sum_{\substack{1\leq r \leq D \\ 1\leq s \leq n}} Y_{\alpha_r}^*\otimes H_{\alpha_{j_s}}\otimes p_{r,j_s}(\mbx). 
\end{equation}
It follows that  \( \varphi_\mbx^{(1)}=
\sum_{\substack{1\leq r \leq D \\ 1\leq s \leq n}} Y_{\alpha_r}^*\otimes H_{\alpha_{j_s}}\otimes p^{(1)}_{r,j_s}(\mbx)\), and \eqref{eq:varphi1} then  yields
\[
\varphi_\mbx^{(1)}( Y_{\alpha_{j_0}}) =
\sum_{s}\, H_{\alpha_s}\otimes p_{j_0,s}^{(1)}(\mbx) = H_{\alpha_{j_0}}\otimes x_{j_0}.
\]
In other words,  \( \varphi_\mbx^{(1)} = \sum_{r=1}^{n} Y^*_{\alpha_{j_r}}\otimes H_{\alpha_{j_r}}\otimes x_{j_r} \), and we conclude that 
\[
q_J(\mbx) \equiv \det \left( p_{j_r,s}^{(1)}(\mbx)  \right) \  = \ x_{j_1}\cdots x_{j_n} = \mbx^J\ \left( \operatorname{mod} \frakm^{n+1} \right),
\]
and hence \( \mbx^J\) lies in \( \cali_\calb\). As a consequence, it follows  that the ideal 
\begin{equation}
\label{eq:ideal_interm}
\bigl< \mbx^J \mid J \in \cals_n(\frg) \bigl>
 \end{equation}
is contained in \(\cali_\calb\).

The final step is to show that the collection \( \{ F_J := \{ \Omega_{\alpha_{j_1}}, \ldots, \Omega_{\alpha_{j_n}} \mid J \in \cals_n(\frg) \}\) consists of  the minimal non-faces of the simplicial complex \( \Delta(\frg)\) on the vertex set \( \scrA(\frg) = \{ \Omega_{\alpha_1}, \ldots, \Omega_{\alpha_D} \} \); see \eqref{eq:CoxArrang}. 

Let \( \{ \alpha_{r_1}, \ldots, \alpha_{r_k} \} \subset R^+\) be a collection of \( k \) distinct positive roots such that 
\( \Omega_{\alpha_{r_1}} \cap \cdots \cap \Omega_{\alpha_{r_k}} \neq \hat \bone := \cap_{\alpha \in R^+} \Omega_\alpha \). Under these conditions, there exists a flat \( \pi \in L(\scrA(\frg))^\circ  \) such that \( \{ \Omega_{\alpha_{r_1}} , \ldots , \Omega_{\alpha_{r_k}}  \} \subset \pi \). In particular, this collection is a face of the complex \( \Delta (\scrA(\frg)) \). Note that this always occurs when \( 1 \leq k < n \) or when \( k = n \) and 
\(\{ \alpha_{r_1}, \ldots, \alpha_{r_n} \} \) is a linearly dependent subset of \( \frh^*_\bbr\). It follows that a non-face needs to have at least \( n \) elements and, in the minimal case, the corresponding collection 
\( \{\alpha_{j_1}, \ldots, \alpha_{j_n} \}\) must be linearly independent. In other words, the element \( J = (j_1, \ldots, j_n) \) lies in \( \cals_n(\frg) \). 

The arguments above can be summarized by saying that the subsets of \( \scrA(\frg) \) which are \emph{minimal non-faces} are in one-to-one correspondence with the set \( \cals_n(\frg) \).  The completes the proof of Theorem \ref{thm:inclusion}.   \qed

\subsection{Proof of Theorem \ref{thm:main}}     
\label{subsec:proof}
Since \( \mu \colon \bbx \to G \) is a morphism of locally finite type between regular schemes, it is a local complete intersection morphism; see \cite[\href{https://stacks.math.columbia.edu/tag/068E}{Lemma 37.62.11}]{stacks-project}. 
On the other hand,  in order to establish the flatness of  \( \mu\), 
it suffices to show\textemdash by Hironaka's criterion\textemdash that its fibers   are equidimensional; see \cite[Cor. 14.130]{MR4225278}.
Furthermore, using faithfully flat descent, we are reduced  to show that the base extension map \( \mu_{\bar\Bbbk} \colon \bbx_{\bar\Bbbk} \to G_{\bar\Bbbk}\) is equidimensional. 

Recall, from \S \ref{subsec:critical}, that \( \scrC \subset \bbx_{\bar\Bbbk}\) is the closed subscheme   of critical points of \(\mu_{\bar\Bbbk}\), and consider the following diagram
\[
\begin{tikzcd}
\bbx_{\bar\Bbbk}  \arrow[d, "\mu_{\bar\Bbbk}"'] \arrow[rr, "\pi_{12}"] & & \bbl_{\bar\Bbbk}\times \bbu_{\bar\Bbbk} \arrow[r,"\rho_2"]& \bbu_{\bar\Bbbk} \arrow[d,"\varphi_{\bar\Bbbk}"] \\
G_{\bar\Bbbk}  & &  & \HomSch{\frl_{\bar\Bbbk}}{\frh_{\bar\Bbbk}},
\end{tikzcd}
\]
where \( \pi_{12} \) is the projection onto the first two factors of 
\( \bbx_{\bar\Bbbk} = \bbl_{\bar\Bbbk}\times \bbu_{\bar\Bbbk}\times \bbl_{\bar\Bbbk}\times \bbu_{\bar\Bbbk}\), and 
\( \rho_2 \) is the projection onto the second factor.

Given a closed point \( p_0 = (x_0, u_0, y_0, v_0) \in \bbx_{\bar\Bbbk} \), it follows from Proposition~\ref{prop:varphi3} that if \( p_0\) lies in \( \scrC \), then
\( u_0 \) lies in \( \calb:= \geo{\varphi}^{-1}(D_{\frl, \fru})\). In particular, one concludes that 
\(
\pi_{12}\left( \scrC\right) \ \subseteq  \rho_2^{-1}\left( \calb \right).
\)
Therefore, 
\begin{equation}
\label{eq:incl}
 \dim \pi_{12}\left( \scrC\right) \ \leq  \dim \rho_2^{-1}\left( \calb \right) = 2D - \operatorname{codim}\left( \calb \right) 
\end{equation}
On the other hand, it follows from Theorem \ref{thm:inclusion} that  
\begin{equation}
\label{eq:codimB}
\operatorname{codim}\left( \calb\right) \geq
 \operatorname{codim}\left( I_{\Delta(\frg)} \right)\geq n = \operatorname{rank }(\frg),
 \end{equation}
 where the last inequality follows directly from Theorem \ref{thm:SRdim}. Therefore,
 \begin{equation}
 \label{eq:final_ineq}
  \dim \pi_{12}\left( \scrC\right) \leq 2D -n = \dim \bbx_{\bar\Bbbk} - \dim G_{\bar\Bbbk} .
\end{equation}

Now, let \( g_0 \) be a closed point in \( G_{\bar\Bbbk}\), and let \( Z \) be an irreducible component of the fiber \( 
 \mu_{\bar\Bbbk}^{-1}(g_0)\). There are two possibilities:
\smallskip

\noindent{\bf 1)}\ \( \boxed{Z \subset \scrC} \) -- 
It follows from Proposition \ref{prop:iso} that the restriction of \( \pi_{12} \) to \( Z\) is a closed embedding. Since
\(\pi_{12}( Z ) \subset \pi_{12}(\scrC)\), one concludes that 
\[ \dim (Z) = \dim( \pi_{12}(Z)) \leq \dim (\pi_{12}(\scrC)) \leq
\dim \bbx_{\bar\Bbbk} - \dim G_{\bar\Bbbk} ;\] see \eqref{eq:final_ineq}. On the other hand,  since \( Z \) is a component of  \( \mu_{\bar\Bbbk}^{-1}(g_0) \) it follows that \( \dim Z \geq \dim \bbx_{\bar\Bbbk} - \dim G_{\bar\Bbbk}  \), and hence equality holds. 
\smallskip

\noindent{\bf 2)}\ \( \boxed{Z \cap \bbx_{\bar\Bbbk}^\circ  \neq \emptyset  }\) -- Recall that \( \bbx_{\bar\Bbbk}^\circ \) is the open   set where \( \mu_{\bar\Bbbk} \) is smooth. It follows that \( Z \cap \bbx_{\bar\Bbbk}^\circ \) is an open dense subset of \( Z \) contained in a fiber of \( \mu_{\bar\Bbbk|  \bbx_{\bar\Bbbk}^\circ}\). Therefore, 
\( \dim Z = \dim (Z\cap \bbx_{\bar\Bbbk}^\circ) = \dim \bbx_{\bar\Bbbk}  -  \dim G_{\bar\Bbbk} \).

This concludes the proof of Theorem \ref{thm:main}.    \qed



\section{Additional results in the \( \scrA_n\) case and open problems}
\label{sec:An}

The main motivation behind this work was the case \( G= SL_{n+1} \), due to its relation to certain constructions in algebraic \( K\)-theory. As a result, this case became the playground where many fun calculations took place and\textemdash at the same time\textemdash the cause of much frustration, due to the fast growth in complexity, for as soon as \( n \geq 4 \) computer packages such as \emph{Macaulay2} could not handle the relevant ideals. 

In this final section we present  additional results  encountered in the \( \scrA_n\). More precisely, we first explain why the fibers of \lulu\ for \(SL_{n+1} \) are global complete intersections. Then we prove that the ideals \( \cali_\calb \) and \( I_{\Delta(\frs\frl_{n+1})} \) coincide, obtaining interesting interpretations for the minimal number of generators for \( \cali_\calb\). Along the way, we
present some open questions in the general case.

%
%

\subsection{Fibers in the \( SL_{n+1} \) case}

Let \( \bbl_n, \bbu_n \subset SL_{n+1} \) denote the unipotent subgroups consisting of the lower (resp. upper) triangular matrices with \( 1\)s along the diagonal, and let \( \bbl\bbu_n \text{ and } \bbu\bbl_n  \) respectively denote the images in \( SL_{n+1} \) of the regular embeddings  
\(
\begin{tikzcd}
\bbl_n\times \bbu_n \arrow[r, hook, "\phi_-"]  &  SL_{n+1}  & \bbu_n \times \bbl_n \arrow[l, hook',"\phi_+"'], 
\end{tikzcd}
\)
 as in \eqref{eq:embed}. 

To describe the ideal defining \( \bbu\bbl_n\), consider   an \( (n+1)\times (n+1) \) matrix of variables  \( \mby = (y_{ij}) \),  and denote by \( \Bbbk[\mby] \) the corresponding polynomial ring, so that 
\( \spec{\Bbbk[\mby] } =: \mat{(n+1)}{(n+1)}  \cong \bba^{(n+1)^2}.\)
Given \( 1 \leq r \leq n+1\), let \( [r]^\complement := [n+1] - [r] \) denote the complement of \( [r]\) in \( [n+1] \), and adopt the convention that \( [0]^\complement = [n+1] \). 
Let
 \( \det_{[r]^\complement, [r]^\complement}(\mby) \in \Bbbk[\mby] \)  be the  minor determinant of the   submatrix of \( \mby \) with columns and rows given by \( [r]^\complement \), and define
\[
\delta_r(\mby) :=  \text{det}_{[r]^\complement, [r]^\complement}(\mby)  - 1 \ \in\  \Bbbk[\mby], \ \  \ r=0, \ldots, n.
\]
It is well-known that \( \bbl\bbu_n\) can be described as follows.

\begin{fact}
\label{fact:LU}
The closed subscheme \( \bbu\bbl_n\) of \( \mat{(n+1)}{(n+1)} \) is given by the ideal \( \scrI_n\) generated by the regular sequence 
\(\{  \delta_0(\mby), \ldots, \delta_{n}(\mby)\}  \ \subset  \Bbbk[\mby]. \)
\end{fact}

\begin{notation}
\label{not:x-u}
Consider the  sets of variables  \( \mbu = (u_{ij} \mid 1 \leq i<j \leq n+1 ) \) and \( \mbx = ( x_{rs} \mid 1 \leq s < r \leq n+1) \),  and  use the same notation \( \mbu= (u_{ij} ) = \left( \begin{smallmatrix}
0 &  u_{ij} \\
 0   & 0 
\end{smallmatrix}\right) \) 
and 
\( \mbx = (x_{r,s} )  = 
\left( \begin{smallmatrix}
0 &  0 \\
x_{rs}   & 0 
\end{smallmatrix}\right)
\)
to denote  the corresponding upper and lower triangular nilpotent matrices. 
For \( k \neq \ell \),    let \( \mbe_{k\ell}  \in \mat{n+1}{n+1}\) be the matrix with \( 1 \) in the \( (k,\ell)\)-entry and \( 0 \) elsewhere. In this way,  the  isomorphisms
\begin{equation}
\label{eq:Psis}
 \Psi_+ \colon \spec{\Bbbk[\mbu]}  \xrightarrow{\ \cong\ } \bbu_n \quad \text{ and } \quad 
  \Psi_- \colon \spec{\Bbbk[\mbx]}  \xrightarrow{\ \cong\ } \bbl_n,
  \end{equation}
induced by the unipotent elements \( \lambda_{k\ell} \colon t \mapsto I + t \mbe_{k\ell} \),  as in \eqref{eq:iso}, are  given precisely by the products
\begin{align*}
\Psi_+(\mbu) & = \prod_{ i<j }^{\text{rev-lex}}\ \lambda_{ij}(u_{ij}) \  =\   I + \mbu 
\quad \text{ and }\quad  \Psi_-(\mbx)  = \prod_{ s<r  }^{\text{lex}}\ \lambda_{rs}(x_{rs}) \  =\   I + \mbx,
\end{align*}
taken in the reverse-lexicographic and lexicographic orders, respectively. 
\end{notation}

\begin{proposition}
\label{prop:CINT}
The fibers of \( \mu \colon \bbl_n\times \bbu_n\times \bbl_
n \bbu_n \to SL_{n+1} \) are global complete intersections. 
\end{proposition}
\begin{proof}
First identify \( \bbl_n\times \bbu_n \) with \( \spec{\Bbbk[\mbu,\mbx]} \),  using the isomorphisms \( \mathbb{\lambda}_+, \mathbb{\lambda}_-\)  above, and 
 consider a closed point \( g_0 \in SL_{n+1} \). We may assume, with no loss of generality, that \( \Bbbk(g_0 ) = \Bbbk \). 
Recall that we showed in Proposition~\ref{prop:iso} that the fiber \( \mu^{-1}(g_0) \) is isomorphic to the subscheme \( \scrV_{g_0} := \phi_-^{-1}(g_0\bbu\bbl_n) \subset \bbl_n \times \bbu_n \). Therefore, for any extension 
 \( \Bbbk  \subset K \), a  \( K\)-valued point \( (\mba, \mbb) \in  \bbl_n(K) \times \bbu_n(K) \) lies in \( \scrV_{g_0}(K) \)  if and only if 
\( g_0^{-1} \mba \mbb \) lies in \( \bbu\bbl_n(K) \).  

Now, define \( f_j^{g_0}(\mbu, \mbx) := \det_{[j]^\complement,[j]^\complement}\left( g_0^{-1} (I+ \mbx)(I+\mbu) \right) - 1 \ \in \ K[\mbu,\mbx] \), for \( j = 1, \ldots, n\) . Note that \( f_0^{g_0}(\mbu, \mbx) = 0\) since all matrices involved are in \( SL_{n+1}\).
Therefore, it follows from Fact \ref{fact:LU} that the ideal of \( \scrV_{g_0} \) in \( K[\mbu, \mbx] \) is generated by 
\( \{  f_1^{g_0}(\mbu, \mbx), \ldots,  f_n^{g_0}(\mbu, \mbx)\} \).  

On the other hand, it is shown in  Theorem \ref{thm:main} that
\( \dim \mu^{-1}(g_0) = 4 \binom{n+1}{2} +1 - (n+1)^2   = n^2\).
Therefore, the codimension of \( \scrV_{g_0} \) in \( \bbl_n \times \bbu_n \) is equal \( 2\binom{n+1}{n} - n^2 = n \), thus proving the proposition. 
\end{proof}

\subsubsection{Open questions}
We have not  been able to answer the following questions.
\begin{enumerate}[\textbf{Q1)}]
\item {Does  \( \{ f_1^{g_0}(\mbu, \mbx), \ldots, f_n^{g_0}(\mbu, \mbx) \} \) become a \emph{regular sequence} after a suitable reordering? It suffices to consider \( g_0 \in N_\bbt ( SL_{n+1} ) \), where 
\(N_\bbt (SL_{n+1})\) is the stabilizer of a maximal torus.}  We can check this with Macaulay2  for  \( n\leq 4\),  for a few choices of \( g_0\).
\end{enumerate}
\begin{enumerate}[\textbf{Q2)}]
\item Are the fibers \( \mu^{-1}(g_0) \) global complete intersections in the  \( \scrB_n, \scrC_n\) and \( \scrD_n\) cases?
\end{enumerate}

\subsection{On the monomial ideals for \( \frs\frl_{n+1}\)}
Once again, we can say a bit more in the case \( \scrA_n\). See Notation \ref{not:x-u} and  contrast the following result with Theorem \ref{thm:inclusion}.

\begin{theorem}
\label{prop:mon-An}
For the semi-simple Lie algebra  \( \frs\frl_{n +1} \),  the following equality   holds:
 \[
  I_{\Delta(\frs\frl_{n+1})} \ =  \  \cali_\calb  . 
 \]
\end{theorem}
Before proving the result, we introduce some auxiliary notation.  Write \( U:= I + \mbu\)\  and\ \( V := (I+\mbu)^{-1}.\)\ \  If \( 1\leq s< r \leq n+1\), then \( \Ad{\Psi_+(\mbu)}{\mbe_{rs}}  = U \, \mbe_{rs} \, V \).
%
Now, let  \( \{ \mbv_k \mid k = 1, \ldots, n\} \) be our fixed basis of  \( \frh = Lie (\bbt)\), with \( \mbv_k:= \mbe_{kk}-\mbe_{n+1, n+1}\). Then we can  write the expression  \eqref{eq:varphi} as 
\begin{align}
\label{eq:explicit}
\varphi_{ \mbu}\   & =\sum_{1\leq i < j \leq n+1} \sum_{k=1}^n  \mbe_{ji}^*\otimes \mbv_k \otimes (  U\mbe_{ji}V)_{kk}
\\ &
=   \sum_{1\leq i < j \leq n+1} \sum_{k=1}^n  \mbe_{ji}^*\otimes \mbv_k \otimes (  V_{ik} U_{kj} ). \notag
\end{align}
Note that \( V = (I+\mbu)^{-1} = \sum_{\ell=0}^n (-1)^\ell \mbu^\ell,\) and hence 
\( V_{ij} = \sum_{\ell=0}^n \ (-1)^\ell  \mbu_{ij}^\ell  \), with
\begin{equation}
\label{eq:uij}
 \mbu_{ij}^\ell \ = \ \sum_{i< a_1< \cdots< a_{\ell-1}< j } \, u_{i,a_1} u_{a_1, a_2} \cdots u_{a_{\ell-1}, j}, \quad \text{for} \quad  \ell \geq 1.
\end{equation}

\begin{notation}
\label{not:indices}
We identify the roots of \( \frs\frl_{n+1} \) with the vectors \( \mbe_j - \mbe_i\) in \( \bbr^{n+1} \), so that  we can think of the following sequence as a set of \( n \) distinct positive roots.  
\begin{equation}
I = \left\{ (i_1, j_1)  \prec \cdots \prec (i_n, j_n) \right\}, \  \text{ with }  1 \leq i_k < j_k \leq n+1 \ \text{ for all }  k = 1, \ldots, n. 
\end{equation}
Given such an \( I \) define:
\begin{enumerate}[a)]
\item \label{it:1} \( \mbm_I(u) := u_{i_1,j_1} u_{i_2, j_2} \cdots u_{i_n, j_n}\), a monomial  in \( \Bbbk[\mbu] \);
\item \label{it:2}  \( \cals(I) := \{ \sigma \in S_n \mid i_k \leq \sigma(k) \leq j_k \} \), a (possibly empty) set of permutations of \( [n]\);
\item \label{it:3}  Given \( \sigma \in \cals(I) \), denote:
\begin{align*}
 B(\sigma) &  := \{ k \in [n] \mid i_k = \sigma(k) < j_k \} \\
 M(\sigma) & := \{ k \in [n] \mid i_k < \sigma(k) < j_k \} \\
 T(\sigma) &  := \{ k \in [n] \mid i_k < \sigma(k) =  j_k \},
\end{align*}
so that \( [n] = B(\sigma)\ \amalg\ M(\sigma)\ \amalg\ T(\sigma) \). 
\item  \label{it:4} We say that \( I \) is \emph{spanning} if \( [n+1] = \{ i_1, \ldots, i_n\} \cup \{ j_1, \ldots, j_n\} \).  This terminology stems from identifying   \( I \) with a set of \emph{edges} of \( \calk_{n+1} \), the complete graph on \( [n+1] \). Therefore, the set \( I \) is spanning iff the corresponding set of edges forms a spanning tree in \( \calk_{n+1} \). 
\item \label{it:5}  \( B_I := \{ \mbe_{j_k} - \mbe_{i_k} \mid k = 1, \ldots, n \} \subset \bbr^{n+1} \), the positive  roots associated to \( I\); 
\item  \label{it:6}  Denote the collection of all such spanning sets by \( \scrG_n \).
\end{enumerate}
\end{notation}
\begin{remark}
\label{rem:facts}
The following statements follow directly from the definitions.
\begin{enumerate}[i.]
\item \label{it:i} If \( I \) is not spanning, then \( \cals(I) = \emptyset \).
\item \label{it:ii} The set \( I \) lies in \( \scrG_n\) iff \( B_I \) is a linearly independent subset of \( \bbr^{n+1} \). 
\item \label{it:iii} The collection \( \left\{ \mbm_I(u) \in \Bbbk[\mbu] \mid I \in \scrG_n \right\} \) is the minimal set of monomial generators of \( I_{\Delta(\frs\frl_{n+1})} \). See  Definition  \ref{def:nerve}.
\end{enumerate}
\end{remark}

\begin{proof}[Proof of  Theorem \ref{prop:mon-An}]
Given \( I = \left\{ (i_1, j_1)  \prec \cdots \prec (i_n, j_n) \right\} \in \scrG_n\), define
\[
 q_I(u) = \det \left(  U_{i_r, k} V_{k, j_r}   \right)_{1\leq r, k \leq n} 
\]
and recall that,  in the \( \frs\frl_{n+1}\) case, the ideal \( \cali_\calb \subset \Bbbk[\mbu] \)  is generated by the polynomials  \(  q_I(u) \) above, as explained in \S\ref{subsec:gens}.

Write
\(
q_I(u) =
\sum_{\sigma \in S_{n}}\varepsilon(\sigma)\,  \textsc{p}_{I,\sigma}(u),
\)
where \( \textsc{p}_{I,\sigma}(u) :=  \prod_{k=1}^n  V_{i_k, \sigma(k)} U_{\sigma(k), j_k}  \)
and observe that, since \( U, V \)  are upper triangular matrices, the polynomial \( \textsc{p}_{I,\sigma}(u) \)  is equal to \( 0 \) unless the permutation \( \sigma  \) lies in \( \cals(I)\).  Therefore, we  can write
\begin{equation}
\label{eq:qI1}
q_I(u) =
\sum_{\sigma \in \cals(I)} \varepsilon(\sigma)\,   \textsc{p}_{I,\sigma}(u).
\end{equation}

Using the partition \( [n] =  B(\sigma)\amalg M(\sigma) \amalg T(\sigma) \), from Notation \ref{not:indices}.\ref{it:4},  we can write
\begin{align*}
& \textsc{p}_{I,\sigma}(u) = \\ 
& =
\left(  \prod_{k \in B(\sigma)} V_{i_k, \sigma(k)} U_{\sigma(k), j_k}\right)
\!\!\! \left(  \prod_{r \in M(\sigma)}  V_{i_r, \sigma(r)}   U_{\sigma(r), j_r}\right) 
\!\!\! \left(  \prod_{s \in T(\sigma)}  V_{i_s, \sigma(s)}  U_{\sigma(s), j_s}\right) \\
& =
\left(  \prod_{k \in B(\sigma)}          u_{i_k, j_k}\right)
\!\!\! \left(  \prod_{r \in M(\sigma)}  V_{i_r, \sigma(r)}   u_{\sigma(r), j_r}\right) 
\!\!\! \left(  \prod_{s \in T(\sigma)}   V_{i_s, j_s} \right),
\end{align*}
since \( U_{ij} = \delta_{ij} + u_{ij} \) and \( V_{ii} = 1 \).

For each \(r \in M(\sigma) \), choose a sequence of integers
\begin{equation}
\label{eq:Ar}
A^r = \{ a^r_0 = i_r < a^r_1 < \cdots < a^r_{\ell_r } < \sigma(r) < j_r \},
\end{equation}
and also denote \( a^r_{\ell_r+1} = \sigma(r) \) and \( a^r_{\ell_r+2} = j_r \). Necessarily, one needs to have \( \ell_{r} +2 < j_r - i_r \). 
Similarly, for each \( s \in T(\sigma) \), choose a sequence of integers
\begin{equation}
\label{eq:Cs}
C^s = \{ c^s_0 = i_s < c^s_1 < \cdots < c^s_{t_s -1} <  c^s_{t_s} =  j_s \},
\end{equation}
where one needs to have \( t_s < j_s - i_s \), as well.

It follows from \eqref{eq:uij} and the considerations above that a typical monomial \( \eta_\star(u) \) appearing as a component of \( \textsc{p}_{I, \sigma}(u) \) has the form
\begin{equation}
\label{eq:nstar}
\eta_\star(u) = 
\hat\mbm_{B(\sigma)}(u) \cdot
\prod_{r \in M(\sigma)} \mbn_{A^r}(u) \cdot
\prod_{s \in T(\sigma)} \mbm_{C^s}(u),
\end{equation}
where
\begin{align}
\hat\mbm_{B(\sigma)}(u) & = \prod_{k \in B(\sigma)} u_{i_k,j_k} \\
 \mbn_{A^r}(u) & =  
\left(  \prod_{p_r=1}^{\ell_r} u_{a^r_{p_r}, a^r_{p_r-1}}\right) \cdot u_{\sigma(r), j_r}  \\
\mbm_{C^s}(u) & = \prod_{q_s=1}^{t_s} u_{c^s_{q_s}, c^s_{q_s-1}} 
\end{align}

\begin{claim}
\label{claim:div}
Given collections \( \{ A^r \mid r \in B(\sigma) \} \) and \( \{ C^s \mid s \in T(\sigma) \} \), as above, one can find \( J \in \scrG_n \) such that the monomial \( \mbm_J(u) \) divides the monomial \( \eta_\star(u) \).
\end{claim}
\begin{proof}[Proof of Claim] 
By hypothesis, the set \( I \) lies in \( \scrG_n\) and by Remark~\ref{rem:facts} one knows that \( \underset{k \in [n]}{\bigwedge }  ( \mbe_{j_k} - \mbe_{i_k} ) \neq 0 \ \in \ \bigwedge^n\bbr^{n+1}.   \) Therefore 
\begin{align}
\label{eq:decomp}
0\neq & \underset{k \in [n]} {\bigwedge} ( \mbe_{j_k} - \mbe_{i_k} ) \\
 & = 
\pm\left\{ \underset{k \in B(\sigma)}{\bigwedge} ( \mbe_{j_k} - \mbe_{i_k} ) \right\} \wedge
\left\{ \underset{r \in M(\sigma)}{\bigwedge} ( \mbe_{j_r} - \mbe_{i_r} ) \right\}\wedge
\left\{ \underset{s \in T(\sigma)}{\bigwedge} ( \mbe_{j_s} - \mbe_{i_s} ) \right\}.   \notag
\end{align}

The middle term in the expression above can be written as
\begin{multline*}
\underset{r \in M(\sigma)}{\bigwedge} \left( \mbe_{j_r} - \mbe_{i_r} \right) 
 =
\underset{r \in M(\sigma)}{\bigwedge} \left\{ \left( \mbe_{j_r} - \mbe_{\sigma(r)} \right) + \left( \mbe_{\sigma(r) } -\mbe_{i_r} \right) \right\} \\
 =
\sum_{D\amalg F = M(\sigma) }\ \varepsilon(D,F)
\left\{\underset{d \in D}{\bigwedge}  \left(\mbe_{j_d} - \mbe_{\sigma(d)}\right) \right\} \wedge \left\{\underset{f \in F}{\bigwedge}\left(\mbe_{\sigma(f) } -\mbe_{i_f}\right) \right\},
\end{multline*}
where \( \varepsilon(D,F) \) is the sign of the shuffle  in the permutation group of \( M(\sigma)\) determined by 
\(D\) and \(F\). 
Hence, it follows from \eqref{eq:decomp} that one can find a partition \( M(\sigma) = D \amalg F \) (with \(D\) or \( F \) possibly empty) such that
\begin{multline}
\label{eq:partition2}
\left\{ \underset{k \in B(\sigma)}{\bigwedge} ( \mbe_{j_k} - \mbe_{i_k} ) \right\} \wedge
\left\{ \underset{d\in D}{\bigwedge} ( \mbe_{j_d} - \mbe_{\sigma(d)} ) \right\}\wedge \\
\wedge\left\{ \underset{f\in F}{\bigwedge} ( \mbe_{\sigma(f)} - \mbe_{i_f} ) \right\}\wedge
\left\{ \underset{s \in T(\sigma)}{\bigwedge} ( \mbe_{j_s} - \mbe_{i_s} ) \right\} \neq 0.
\end{multline}
%
%

Now, consider choices \(\{ A^r \mid r \in M(\sigma) \} \) and \( \{ C^s \mid s \in T(\sigma) \} \) as in \eqref{eq:Ar} and \eqref{eq:Cs}, respectively, giving rise to the ``typical monomial'' \( \eta_\star(u) \); see \eqref{eq:nstar}.  Note that one can write
\begin{align*}
\underset{f \in F}{\bigwedge}& \left( \mbe_{\sigma(f)} - \mbe_{i_f} \right) =
\underset{f \in F}{\bigwedge} \left\{ \sum_{p_f=1}^{\ell_f} \left( \mbe_{a^f_{p_f}} - \mbe_{a^f_{p_f-1}}\right) \right\},
\end{align*}
and, similarly, 
\[
\underset{s \in T(\sigma)}{\bigwedge} \left( \mbe_{j_s} - \mbe_{i_s} \right) =
\underset{s \in T(\sigma)}{\bigwedge} \left\{ \sum_{q_s=1}^{t_s} \left( \mbe_{c^s_{q_s}} - \mbe_{c^s_{q_s-1}} \right)\right\}.
\]
Then,  it follows from \eqref{eq:partition2} that for each \( f \in F \) one can find \( 1 \leq p^0_f \leq \ell_f  \), and for each \( s \in T(\sigma) \) one can find \( 1 \leq q^0_s \leq t_s \), such that

\begin{multline}
\label{eq:partition3}
\left\{ \underset{k \in B(\sigma)}{\bigwedge} \left( \mbe_{j_k} - \mbe_{i_k} \right) \right\} \wedge
\left\{ \underset{d\in D}{\bigwedge} \left( \mbe_{j_d} - \mbe_{\sigma(d)} \right) \right\}\wedge \\
\wedge \left\{ \underset{f\in F}{\bigwedge} \left( \mbe_{a^f_{p_f^0}} - \mbe_{a^f_{p_f^0-1}} \right) \right\}\wedge
\left\{ \underset{s \in T(\sigma)}{\bigwedge} \left(  \mbe_{a^s_{q_s^0}} - \mbe_{c^s_{q_s^0-1}}  \right) \right\} \neq 0.
\end{multline}

The conclusion above, along with Remark \ref{rem:facts}.\ref{it:ii}, implies that disjoint union
\begin{align*}
J & := \left\{ (i_k, j_k)\right\}_{ k \in B(\sigma) }\ \coprod \
 \left\{ (\sigma(d), j_d)\right\}_{ d \in D} 
\ \coprod \ \ \left\{ (i_f, \sigma(f) ) \right\}_{f \in F }
\ \coprod \  \left\{ (i_s, j_s )\right\}_{s \in T(\sigma) }
\end{align*}
forms a spanning set in the sense of Notation \ref{not:indices}.\ref{it:4}, i.e. \( J \) lies in \( \scrG_n\).

Now, since \( M(\sigma) = D\amalg F \), if follows from \eqref{eq:nstar} that \( \mbm_J(u)  \) divides the monomial \( \eta_\star(u) \). This concludes the proof of Claim \ref{claim:div}.
\end{proof}

Recall that the monomial 
\[
\mbm_J (u) := \prod_{k\in B(\sigma)} u_{i_k, j_k} \cdot \prod_{d \in D} u_{\sigma(d),  j_d } \cdot
\prod_{f\in F} u_{a^f_{p_f^0-1}, a^f_{p_f^0}} \cdot
\prod_{s\in T(\sigma)} u_{c^s_{q_s^0-1}, c^s_{q_s^0}}
\]
lies in \( I_{\Delta(\frs\frl_{n+1})} \), according to Remark \ref{rem:facts}.\ref{it:iii}.  On the other hand, it follows from Claim \ref{claim:div} that all  monomials \( \eta_\star(u) \)  that are components of \( q_I(u) \) lie in \( I_{\Delta(\frs\frl_{n+1})} \). 
The theorem follow from this fact and \S\ref{subsec:gens}.
\end{proof}

\subsubsection{Open questions.} The answer to the following questions would be of interest.

\begin{enumerate}[\textbf{Q3)}]
\item  Does the equality  \( I_{\Delta(\frg)}  = \cali_\calb \) hold in the  \( \scrB_n, \scrC_n\) and \( \scrD_n\) cases?
\end{enumerate}

We have seen that the cardinality of the minimal set of monomial generators in the case \( \scrA_n\), for the ideals \( I_{\Delta(\frg)} = \cali_\calb \), is the number of spanning trees of  \( \calk_{n+1} \), the complete graph on \( [n+1] \).  See
Notation \ref{not:indices}.\ref{it:4}. This number is precisely \( \mba(n) := (n+1)^{n-1} \), and this sequence of numbers has many other combinatorial interpretations in the literature. 

Now, let \( W \) be a finite Coxeter group, with irreducible Dynkin diagram \( D \).  In a letter \cite{Del-Loo74} to E. Looijenga, Pierre Deligne calculated\textemdash among other things\textemdash the number of distinct ways one can express a Coxeter element \( c \in W\)  as a product of  \( n = \operatorname{rank}(D) \) reflections, and he denoted this number by \( F'(D)\). It turns out that,  in the \( \scrA_n\) case, this number is precisely \( \mba(n)= (n+1)^{n-1}\).

\begin{enumerate}[\textbf{Q4)}]
\item  What is the cardinality of the minimal set of monomial generators for \( \cali_\calb\) in the \( \scrB_n, \scrC_n\) and \( \scrD_n\) cases? Same question for \( I_{\Delta(\frg)} \) in these cases.  
\end{enumerate}
\begin{enumerate}[\textbf{Q5)}]
\item Is there a relationship between the answers to the previous question and Deligne's \( F'(D) \) in these cases?
\end{enumerate}

\bibliographystyle{amsalpha}

\providecommand{\bysame}{\leavevmode\hbox to3em{\hrulefill}\thinspace}
\providecommand{\MR}{\relax\ifhmode\unskip\space\fi MR }
\providecommand{\MRhref}[2]{%
  \href{http://www.ams.org/mathscinet-getitem?mr=#1}{#2}
}
\providecommand{\href}[2]{#2}

%

\end{document}